\documentclass[article]{elsarticle}

\usepackage{hyperref}

\usepackage{amsfonts,amssymb,amsmath}
\usepackage{graphicx}
\newcommand{\R}{\mathbb{R}}
\newcommand{\N}{\mathbb{N}}

\newcommand{\dd}[1]{\mathrm{d}#1}

\newtheorem{theorem}{Theorem}[section]
\newtheorem{lemma}{Lemma}[section]
\newtheorem{remark}{Remark}[section]
\usepackage{a4wide,color}
\usepackage{subfigure}

\journal{}

\begin{document}

\begin{frontmatter}

\title{
Fractional time stepping and adjoint based gradient computation in an inverse problem for a fractionally damped wave equation 
}
\tnotetext[mytitlenote]{Supported by FWF under grants P30054 and DOC 78}

\author{Barbara Kaltenbacher and Anna Schlintl}
\address{University of Klagenfurt, Austria}

\begin{abstract}
In this paper we consider the inverse problem of identifying the initial data in a fractionally damped wave equation from time trace measurements on a surface, as relevant in photoacoustic or thermoacoustic tomography. We derive and analyze a time stepping method for the numerical solution of the corresponding forward problem. Moreover, to efficiently obtain reconstructions by minimizing a Tikhonov regularization functional (or alternatively, by computing the MAP estimator in a Bayesian approach), we develop an adjoint based scheme for gradient computation. Numerical reconstructions in two space dimensions illustrate the performance of the devised methods. 
\end{abstract}


\end{frontmatter}


\section{Introduction}
Inverse problems for the classical second order (acoustic) wave equation or its elastic or electromagnetic counterpart are a highly active field of research since they arise in numerous applications ranging from ultrasound imaging via geophysical prospection to microwave tomography.
We will here particularly focus on the problem of reconstructing the initial data from time trace measurements, as it arises, for example, in photoacoustic (PAT) or thermoacoustic tomography (TAT).  
From a regularization point of view, these inverse problems are only mildly ill-posed as long as the wave propagation is lossless, due to the fact that it can then be basically inverted without losing information.
However, as soon as attenuation is involved, the situation is different: Strong damping is known to render the inverse problem even severely ill-posed. Recently, fractional order damping models have been put forward due to their physical relevance in, e.g., ultrasound propagation. The frequency dependence of attenuation follows a power law, consequently fractional order time derivatives govern the damping in the time domain setting, see, e.g., the survey \cite{CaiChenFangHolm_survey2018} and the references therein.

Following up on \cite{kaltenbacherrundell:fracpat}, we consider PAT, taking power law frequency dependent attenuation into account by using the Caputo-Wismer-Kelvin model \cite{Caputo:1967,Wismer:2006}.
This amounts to the inverse problem of identifying $u_0$ in the fractionally damped wave equation 
\begin{equation}
\label{mod_init}
\begin{aligned}
u_{tt} -c^2 \Delta u - b \Delta \partial_t^{\alpha} u &= 0 \mbox{ in } D \times (0,T) \\
u(0) = u_0, \quad u_t(0) &= 0 \mbox{ in } D,
\end{aligned}
\end{equation}
where $u(x,t)$ denotes the acoustic pressure, $c>0$ the speed of sound, and $b>0$ the diffusivity of sound. 
For more details on photoacoustic tomography in the lossless case, we refer, e.g, to the review paper
\cite{KuchmentKunyanski:2011} and the references therein; for the attenuated case, see, e.g., 
\cite{Ammari2012,ElbauScherzerShi:2017,KowarScherzer:2011}.


Here $D$ is either all of $\mathbb{R}^d$, $d\in\{1,2,3\}$ (with some radiation condition imposed on $u$) or $D\subset \mathbb{R}^d$ is a bounded $C^{1,1}$ smooth or convex polygonal domain. In the latter case, boundary conditions are imposed on $u$. If $D$ is large enough, these may just be homogeneous Dirichlet or Neumann conditions; otherwise, in order to avoid spurious reflections on $\partial D$, absorbing boundary conditions or so-called perfectly matched layers should be employed.   

The time differential operator $\partial_t^\alpha$ in this model is the Djrbashian-Caputo fractional time derivative, which is essential in order to allow for prescribing initial conditions.
For details on fractional differentiation and subdiffusion equations, we refer to, e.g.,
\cite{Dzharbashyan:1966t,Djrbashian:1993,MainardiGorenflo:2000,SakamotoYamamoto:2011a,SamkoKilbasMarichev:1993},
see also the tutorial on inverse problems for anomalous diffusion processes
\cite{JinRundell:2015}.
The Caputo-Djrbashian derivative of order $\alpha \in (n-1,n)$ with $n \in \N$ of a function $u:[0,T) \to \R$ is defined by
\[
D_t^{\alpha} u = I^{n-\alpha} [u^{(n)}],
\]
where $u^{(n)}$ denotes the $n$-th integer order derivative and for $\gamma \in (0,1)$ the Abel integral operator $I^{\gamma}$ is defined by
\begin{equation}\label{eqn:Abel}
I^{\gamma}[v](t) = \frac{1}{\Gamma(\gamma)} \int_0^t \frac{v(s)}{(t-s)^{1-\gamma}} \dd{s}.
\end{equation}

The additional observations available in PAT for determining the unknown $u_0$ are measurements of the sound pressure level at distance from the object to be imaged
\begin{equation}\label{eqn:observations}
u= y \quad \mbox{ on }\Sigma\times(0,T)
\end{equation}
with a closed surface $\Sigma=\partial\Omega$ that is entirely contained in the interior of $D$.

In \cite{kaltenbacherrundell:fracpat}, the inverse problem of identifying $u_0$ in \eqref{mod_init} from observations \eqref{eqn:observations} has been studied. Well-posedness of the forward problem, that is, the initial boundary value problem for \eqref{mod_init} with given $u_0$, as well as uniqueness for the inverse problem have been established with the above and also with different damping models, and reconstructions in one space dimension have been provided. 

In this paper, we focus on two key computational aspects of this inverse problem, namely on a time stepping method for the numerical solution of the forward problem and on adjoint based gradient computation for the reconstruction of $u_0$. 
Both tasks require special treatment of the fractional damping term and both aspects are relevant more generally in inverse problems for fractionally damped wave equations, so beyond the PAT/TAT context we are focusing on in our numerical experiments.

\subsection{Reconstruction by Bayesian inference or Tikhonov regularization}\label{sec:bayes}
Linear inverse problems, like the one under consideration here, can be written as operator equations
\begin{equation}
\label{modeleq}
y^{\delta} = \mathcal{G} u_0 + \delta \eta,
\end{equation}
where $\mathcal{G}$ is the (linear) forward operator, mapping between Hilbert spaces $X$ and $Y$, $\delta\ge0$ describes the noise level, $u_0$ is the unknown target. In the inverse problem described above, the forward operator $\mathcal{G}$ is the map that takes $u_0$ to the boundary trace according to \eqref{eqn:observations} of the solution $u$ to \eqref{mod_init}.  Moreover, $y^{\delta}$ is the given noisy data, with $\eta$ a random variable describing the noise. We restrict our considerations to Gaussian noise, i.e. $\delta\eta \sim \mathcal{N}(0, \Gamma_{noi})$. The prior distribution for $u_0$ is chosen to be normal as well, i.e. $u_0 \sim \mathcal{N}(u_0^*,  \Gamma_{pr})$. In that setting, the posterior is also normally distributed with mean and covariance given by
\[
{u_0}_{\alpha}^{\delta} = (\mathcal{G}^*\Gamma_{noi}^{-1}\mathcal{G}+\Gamma_{pr}^{-1})^{-1} \bigl(\mathcal{G}^*\Gamma_{noi}^{-1}y^\delta+\Gamma_{pr}^{-1}u_0^*\bigr)\,,\qquad
\Gamma_{post} = (\mathcal{G}\Gamma_{noi}^{-1}\mathcal{G}^*+\Gamma_{pr}^{-1})^{-1},
\]
see, e.g., \cite[Theorem 2.4]{stuart:2010}.
Here $u_0^*$ is an element of $X$, $\Gamma_{pr}:X\to X$, $\Gamma_{noi}:Y\to Y$ are selfadjoint positive definite operators, and the superscript ${}^*$ denotes the Hilbert space adjoint. 

Setting up the matrices for computing ${u_0}_{\alpha}^{\delta}$ can be avoided by applying a gradient descent method (possibly accelerated by a quasi Newton preconditioner such as limited memory BFGS) to the minimization problem of computing ${u_0}_{\alpha}^{\delta}$ as a maximum a posteriori MAP estimator 
\begin{equation}\label{MAP}
\min_{u_0\in X} \frac12\|\Gamma_{noi}^{-1/2}(\mathcal{G}u_0-y^\delta)\|^2_Y+\frac12\|\Gamma_{pr}^{-1/2}(u_0-u_0^*)\|_X^2
\end{equation}
Clearly, without computing the covariance, this is plain Tikhonov regularization.

To numerically solve the optimization problem of finding the MAP estimator, for example by some quasi Newton method it is essential to compute the gradient of the cost function $J(u_0)=\frac12\|\Gamma_{noi}^{-1/2}(\mathcal{G}u_0-y^\delta)\|^2_Y+\frac12\|\Gamma_{pr}^{-1/2}(u_0-u_0^*)\|_X^2$ in \eqref{MAP}
\begin{equation}\label{gradJ}
J'(u_0)= \mathcal{G}^*\Gamma_{noi}^{-1}(\mathcal{G}u_0-y^\delta)+\Gamma_{pr}^{-1}(u_0-u_0^*)
\end{equation}
For this purpose we will use an adjoint approach which safes computational effort by avoiding the numerical solution of sensitivity equations.
For this purpose, in Section \ref{sec:adjoint}, we derive the adjoint of the forward operator and provide an algorithm for computing the gradient of $J$.

\subsection{Numerical solution of the forward problem}
In each step of some iterative reconstruction scheme, (for example a gradient based one, as indicated above) the initial boundary value problem \eqref{mod_init} for a fractionally damped linear wave equation needs to be solved. 
After space discretization (e.g. with finite elements), \eqref{mod_init} 
becomes a system of ODEs 
\begin{equation} \label{eqn:semidiscr}
M \ddot{u}  + \partial_t^{\alpha} Cu + Ku = f,
\end{equation}
where the dot denotes the time derivative.
The self-adjoint positive (semi-)definite matrices $M$, $K$, and $C$, represent discretizations of the identity ($M$), the negative Laplacian ($K$) and the damping ($C$), which in case of \eqref{mod_init} is identical to $K$, but might as well represent, e.g. a combination of $M$ and $K$ (Rayleigh damping) or a fractional Laplacian (e.g. in the Chen-Holm model \cite[eq. (21)]{ChenHolm:2004}, see also \cite{kaltenbacherrundell:fracpat}.) 
To numerically solve the semidiscrete system \eqref{eqn:semidiscr}, in Section \ref{sec:timestepping} we derive and analyze a time stepping scheme, that is based on the Newmark method, equipped with a quadrature formula for the fractional derivative. 

\medskip

The remainder of this paper is organized as follows.
Section \ref{sec:timestepping} focuses on the numerical solution of the forward problem, in particular, derivation and analysis of a time stepping scheme. We establish a stability estimate and convergence with rates for sufficiently smooth solution. 
In Section \ref{sec:adjoint} we derive the adjoint model that is used in the numerical reconstruction scheme consisting of a gradient based method for computing the MAP estimator.
Finally, Section \ref{sec:reconstructions} provides spatially two-dimensional reconstructions based on an implementation of the previously derived computational tools.

\def\basis{N}
\section{Time stepping scheme}\label{sec:timestepping}
\subsection{Newmark method}
For the time discretization of the problem, we are going to adapt the Newmark scheme (which is also known to be a particular case of the generalized alpha scheme \cite{ChungHulbert1993}) to enable inclusion of the fractional derivative. Space discretization is supposed to be done beforehand by a finite element method 
$u(x)\approx \sum_{i=1}^{n_{el}} u_i \basis_i(x)$ with basis functions $\basis_i$, so that $M$ and $K$ are the mass and stiffness matrices according to 
\[
M_{i,j}=\int_\Omega \basis_i(x)\basis_j(x) \, dx,\qquad 
K_{i,j}=\int_\Omega \nabla\basis_i(x)\cdot\nabla\basis_j(x) \, dx. 
\]
While for \eqref{mod_init}, we simply have $C=K$, the framework of this section actually allows for an  arbitrary positive (semi-)definite damping matrix $C$.

From now on we denote by $u_{n+1}, \dot{u}_{n+1}, \ddot{u}_{n+1}$ the approximations for the function and derivative values $u(t_{n+1}), \dot{u}(t_{n+1}), \ddot{u}(t_{n+1})$. 
The time discretization of the fractional derivative, which is nonlocal in time, more precisely, given by a convolution integral, will be approximated by a discrete convolution $\partial_t^\alpha u(t_n) \approx \sum_{j=0}^{n}b_j^{n} \dot{u}(t_{n-j})$. 
There exists a vast amount of literature on this task, see, e.g., the overview provided in \cite{jin:overview}.
We will here consider two approaches for this purpose, namely an L1 type scheme and a Galerkin scheme in sections \ref{sec:L1} and \ref{sec:Galerkin}, respectively.

The Newmark scheme approximates \eqref{eqn:semidiscr} by the time discrete model
\begin{equation}
\label{Newmark}
M \ddot{u}_{n+1} +  C \sum_{j=0}^{n+1}b_j^{n+1} \dot{u}_{n+1-j} + Ku_{n+1} = f_{n+1}
\end{equation}
where
\[ \begin{aligned}
u_{n+1} &= u_n + \Delta t \dot{u}_n + \frac{\Delta t^2}{2} ((1-2\beta)\ddot{u}_n + 2\beta \ddot{u}_{n+1}) \\
\dot{u}_{n+1} &= \dot{u}_n + \Delta t ((1-\gamma)\ddot{u}_n + \gamma \ddot{u}_{n+1}.
\end{aligned} \]
Since it is an implicit method, we have to solve the following system in each time step
\[ \begin{aligned}
&M^* \ddot{u}_{n+1} = f_{n+1} -  C\left( b_0^{n+1} (\dot{u}_n + (1-\gamma) \Delta t \ddot{u}_n) + \sum_{j=1}^{n+1} b_j^{n+1} \dot{u}_{n+1-j}\right) \\
&\hspace*{5cm}- K (u_n + \Delta t\dot{u}+ \frac{\Delta t^2}{2}(1-\beta)\ddot{u}_n) \\
&\mbox{where }M^* = M + b_0^{n+1}\gamma \Delta t C + \beta \Delta t^2 K.
\end{aligned} \]

To solve this system, we consider two formulations (in analogy to the commonly used formulations in the integer derivative case), an effective mass formulation and an effective stiffness formulation. In the former we solve for the second derivative $\ddot{u}$, while in the latter we solve directly for $u$.

\underline{Effective mass formulation}\\
Predict:
\[ \begin{aligned}
\tilde{u} &= u_n + \Delta t \dot{u}_n + (1-2\beta) \frac{\Delta t^2}{2} \ddot{u}_n \\
\tilde{\dot{u}} &= \dot{u}_n + \Delta t(1-\gamma) \ddot{u}_n
\end{aligned} \]
Solve:
\begin{equation}\label{solve_mass}
\begin{aligned}
&M^* \ddot{u}_{n+1} = f_{n+1} - K \tilde{u} -  C \left(b_0^{n+1} \tilde{\dot{u}} + \sum_{j=1}^{n+1} b_j^{n+1} \dot{u}_{n+1-j}\right) \\
&\mbox{where }M^* = M + b_0^{n+1}\gamma \Delta t C + \beta \Delta t^2 K
\end{aligned}
\end{equation}
Correct:
\begin{equation}\label{corr_mass}
\begin{aligned}
u_{n+1} &= \tilde{u} + \beta \Delta t^2 \ddot{u}_{n+1} \\
\dot{u}_{n+1} &= \tilde{\dot{u}} + \gamma \Delta t \ddot{u}_{n+1} 
= \dot{u}_n + \Delta t((1-\gamma) \ddot{u}_n + \gamma \ddot{u}_{n+1})
\end{aligned}
\end{equation}

\underline{Effective stiffness formulation}\\
Predict:
\[ \begin{aligned}
\tilde{u} &= u_n + \Delta t \dot{u}_n + (1-2\beta) \frac{\Delta t^2}{2} \ddot{u}_n \\
\tilde{\dot{u}} &= \dot{u}_n + \Delta t(1-\gamma) \ddot{u}_n
\end{aligned} \]
Solve:
\[ \begin{aligned}
&K^* u_{n+1} = f_{n+1} - C \left( b_0^{n+1}(\Delta t \tilde{\dot{u}} - \frac{\gamma}{\beta} \tilde{u} ) + \sum_{j=1}^{n+1} b_j^{n+1} \dot{u}_{n+1-j}\right) + M \frac{1}{\Delta t^2 \beta} \tilde{u} \\
&\mbox{where }K^* = K + b_0^{n+1}\frac{\gamma}{\beta} C + \frac{1}{\beta \Delta t^2} M
\end{aligned} \]
Correct:
\[ \begin{aligned}
\dot{u}_{n+1} &= \tilde{\dot{u}} + \frac{\gamma}{\beta \Delta t} (u_{n+1} - \tilde{u}) = \tilde{\dot{u}} + \gamma \Delta t \ddot{u}_{n+1} \\
\ddot{u}_{n+1} &= \frac{u_{n+1}-\tilde{u}}{\beta \Delta t^2}
\end{aligned} \]

\subsection{L1 type scheme} \label{sec:L1}
The fractional derivative is discretized similarly to the L1 scheme \cite{jin:overview} but relying on values of $\dot{u}$ rather than $u$.
\[ \begin{aligned}
&\partial_t^{\alpha} u(t_n) 
= \frac{1}{\Gamma ( 1 - \alpha)} \sum_{j=0}^{n-1} \int_{t_j}^{t_{j+1}} \dot{u}(s) (t_n -s)^{-\alpha} \dd{s} 
 \approx \frac{1}{\Gamma (1- \alpha)} \sum_{j=0}^{n-1} \frac{\dot{u}(t_{j+1}) + \dot{u}(t_j)}{2} \int_{t_j}^{t_{j+1}} (t_n -s)^{-\alpha} \dd{s} \\
&= \frac{1}{2\Gamma (1- \alpha)} \sum_{j=0}^{n-1} \left(\dot{u}(t_{j+1}) + \dot{u}(t_j)\right) \frac{1}{1-\alpha} \left( (t_n - t_j)^{1-\alpha} - (t_n - t_{j+1} )^{1-\alpha} \right) \\
& = \frac{1}{2 \Gamma (2- \alpha)} \left\{ \sum_{j=0}^{n-1} \dot{u}(t_j) \left( (t_n - t_j)^{1-\alpha} - (t_n - t_{j+1})^{1-\alpha} \right) + \sum_{i=1}^n \dot{u}(t_i) \left( (t_n - t_{i-1})^{1-\alpha} - (t_n - t_i)^{1-\alpha} \right) \right\}
\end{aligned} \]
and we obtain for the discretized fractional derivative
\[ \begin{aligned}
\partial_t^{\alpha} u(t_n) \approx \sum_{j=0}^n b_j^n \dot{u}(t_{n-j})
\end{aligned} \]
with
\begin{equation}\label{L1scheme} 
\begin{aligned}
t_j = j \Delta t \quad \mbox{ and } \quad  b_j^n = \frac{\Delta t^{1-\alpha}}{2\Gamma(2-\alpha)}\left\{ \begin{array}{cc}  (n^{1-\alpha} -(n-1)^{1-\alpha})  & j=n \\
1 & j = 0 \\
((j+1)^{1-\alpha} - (j-1)^{1-\alpha})  & j \in \{1,\ldots,n-1\} \end{array} \right.
\end{aligned} 
\end{equation}

\subsubsection{Comparison to an exact solution}

To validate the algorithm and its implementation we construct an example that allows to compute the exact solution. We will rely on the solution representation via an expansion in terms of eigenfunctions $\psi_j$ of $-\Delta$ 
\[
u(x,t) = \sum_{j=1}^{\infty} u_j(t) \psi_j(x),
\]
where for each $j \in \N$, $u_j$ solves the relaxation equation
\begin{equation}\label{relaxationeq}
\ddot{w}_j + c^2 \lambda_j u_j + b \lambda_j \partial_t^{\alpha} w =: \ddot{w} + A \partial_t^{\alpha}w + Bw = 0 
\end{equation}
with initial conditions $w(0) = 1$ and $\dot{w}(0)=0$, and $\lambda_j$ the eigenvalue corresponding to $\psi_j$, see also \cite{kaltenbacherrundell:fracpat}. 
We make use of the explicit solution of \eqref{relaxationeq} computed in  \cite{sakakibara:order12} for the special case $\alpha := \frac{1}{2}$, $A:= \frac{4}{3^{3/4}}$ and $B :=1$. The Laplace transform yields
\begin{equation*}
\hat{w}(s) = \frac{s+\frac{4}{3^{3/4}} s^{-1/2}}{s^2+\frac{4}{3^{3/4}} s^{1/2}+1}.
\end{equation*}
The denominator $\omega(s) := s^2+\frac{4}{3^{3/4}} s^{1/2} + 1$ has two zeros, namely
\begin{equation*}
p_+ = \frac{1}{3}(-\sqrt{3} + 2i \sqrt{6}), \quad p_{-} = \frac{1}{3}(-\sqrt{3} - 2i \sqrt{6}).
\end{equation*}
Thus, we have
\begin{equation*}
\mathcal{R}_{1/2}(t) = 2 Re \left\{ \frac{p_+ + A p_{+}^{-1/2}}{2p_+ + \frac{1}{2}Ap_+^{-1/2} } e^{p_+ t} \right\} = \frac{1}{3} e^{-t/\sqrt{3}} \left( \cos (2 \sqrt{\frac{2}{3}}t) + \frac{1}{\sqrt{2}} \sin(2 \sqrt{\frac{2}{3}}t) \right)
\end{equation*}
and the spectral function
\begin{equation*}
\begin{aligned}
H_{1/2}(r) &= \frac{A \sin( \pi /2)}{\pi} \frac{r^{-1/2}}{(r^2+ Ar^{1/2} \cos(\pi /2) + 1)^2 + (A r^{1/2} \sin ( \pi /2))^2} = \frac{4}{3^{3/4} \pi r^{1/2} ((r^2+1)^2+\frac{16r}{3 \sqrt{3}})}.
\end{aligned}
\end{equation*}
Finally, we obtain the exact solution which is given by
\begin{equation*}
w(t) = \mathcal{R}_{1/2}(t) + \int_0^{\infty} e^{-rt} H_{1/2}(r) \dd{r}.
\end{equation*}

\subsubsection{Numerical experiment with L1 type scheme}
For the experiment we use as initial condition the first eigenfunction on the unit square which along with its eigenvalue is given by
\[
\psi_1(x,y) = \sin( \pi x) \sin(\pi y), \quad \lambda_1 = 2\pi^2,
\]
resulting in the exact solution
\[
u(t) = w(t) \sin ( \pi x ) \sin (\pi y)
\]
with $w$ according to \eqref{relaxationeq}, $w(0)=1$, $\dot{w}(0)=0$
The numerical solution is computed by the above described extended Newmark scheme with the L1 scheme coefficients according to \eqref{L1scheme} with the parameters $c^2 = \frac{1}{\pi^2}$ and $b = \frac{4}{3^{3/4} \pi^2}$ resulting from the choice of $A$ and $B$ respectively. The parameters for the Newmark scheme are chosen to be $\gamma := 0.5$ and $\beta := 0.25$, so optimally with respect to convergence order. The initial values are given by
\[ \begin{aligned}
u_0 &= \sin( \pi x) \sin (\pi y) \\
\dot{u}_0 &= 0 \\
M \ddot{u}_0 &= (f - b C\partial_t^\alpha u - c^2 K u)(t=0) = f(0) - c^2 K u_0.
\end{aligned} \]

The numerical solution is then computed with $T=4$ and $\Delta t = 0.08$, i.e. we have 50 time steps. The following result show the numerical solution in the left column and the exact solution in the right column at times $0, 0.8, 1.6, 2.4, 3.2, 4$. 

\begin{figure}[ht]
\includegraphics[width=0.49\textwidth]{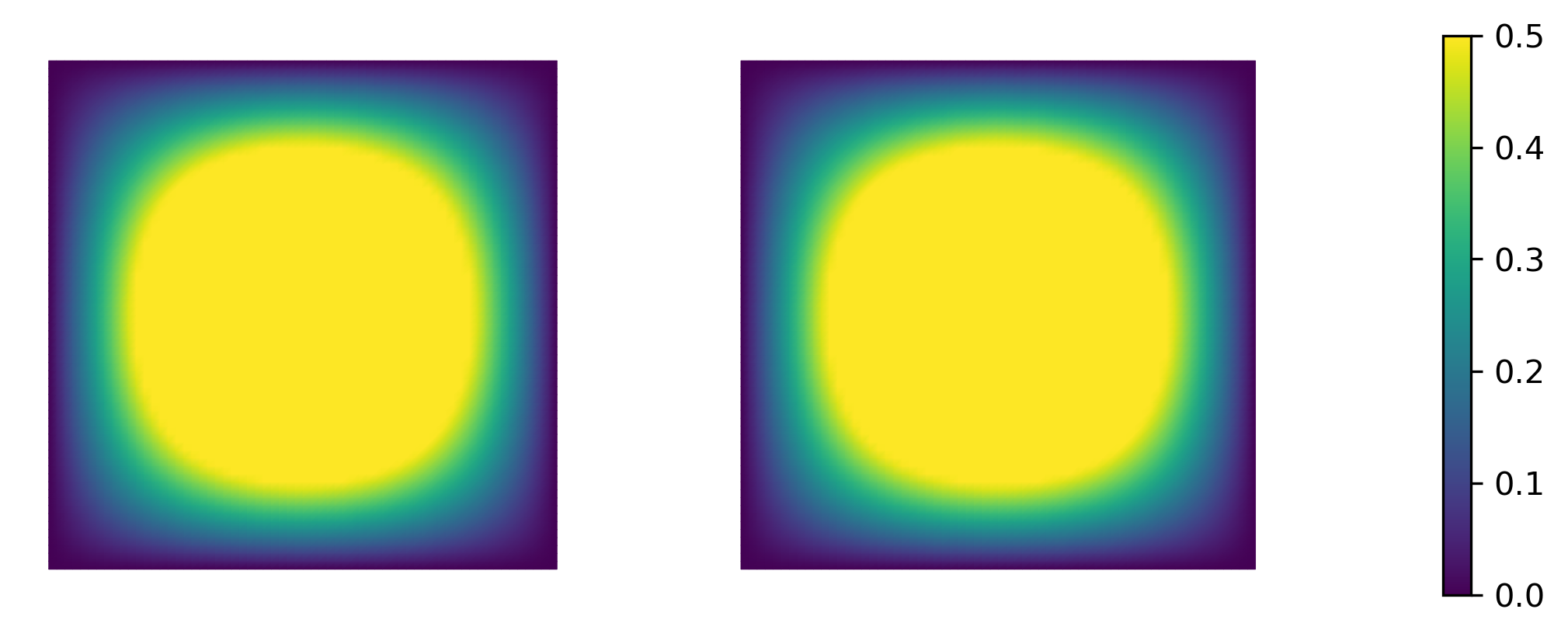}
\includegraphics[width=0.49\textwidth]{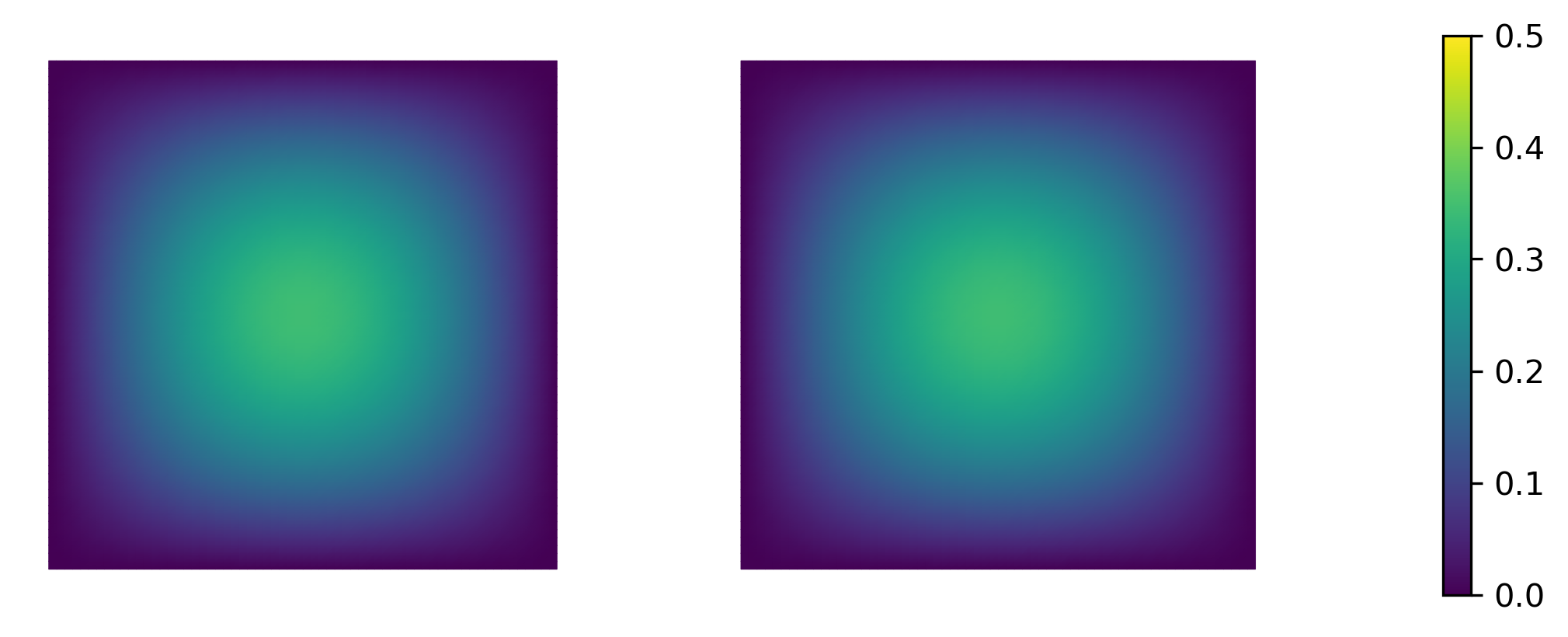}\\
\includegraphics[width=0.49\textwidth]{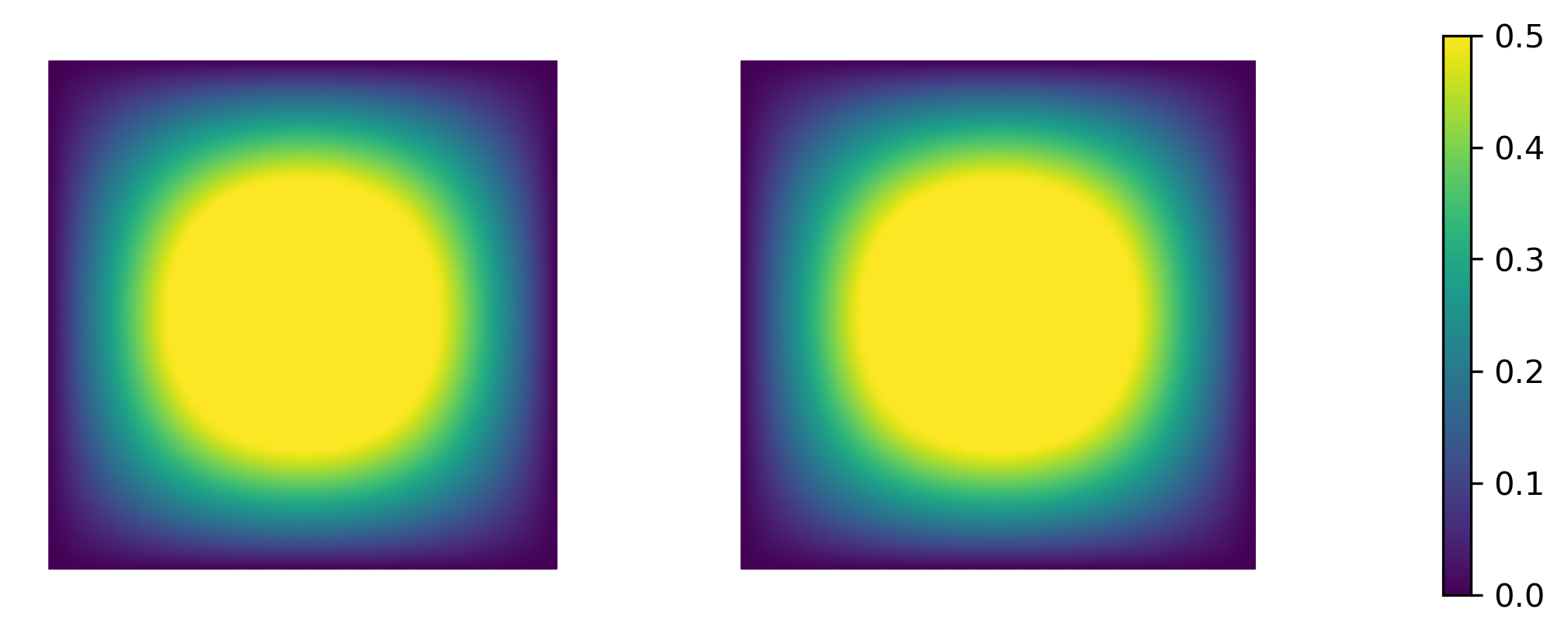}
\includegraphics[width=0.49\textwidth]{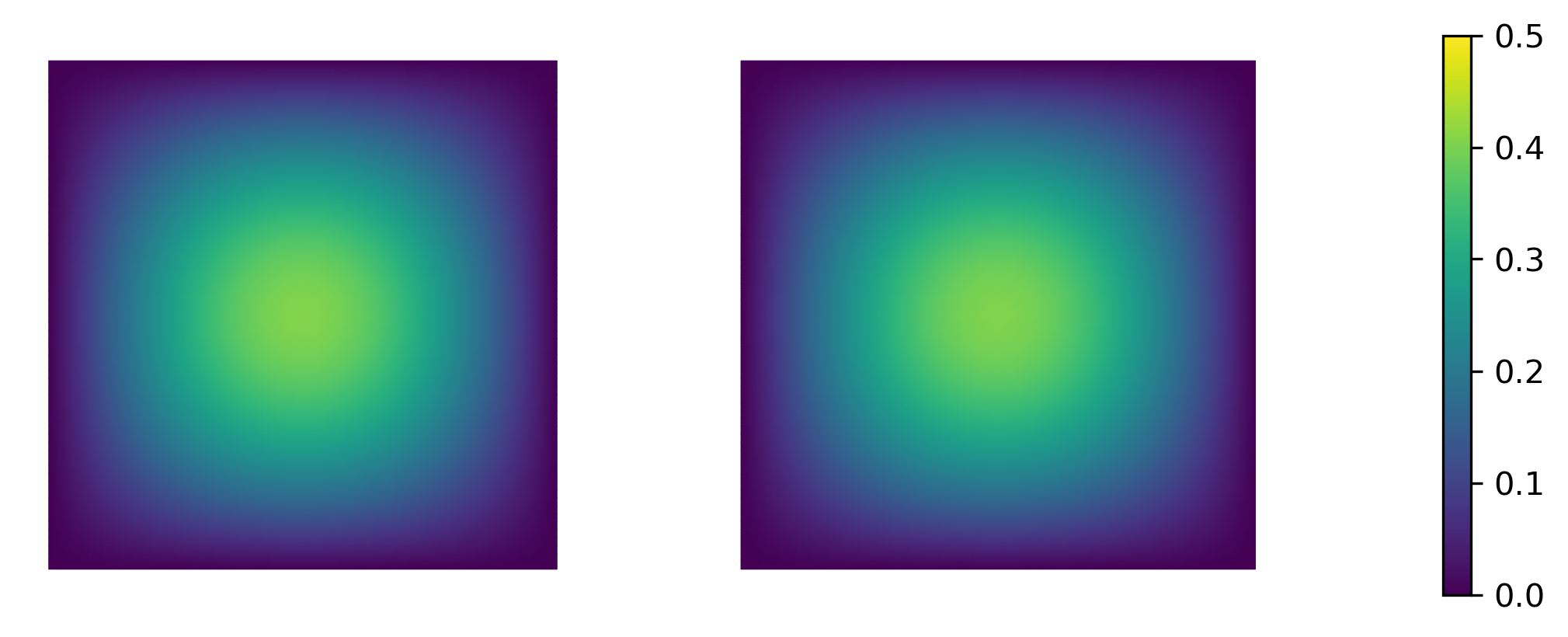}\\
\includegraphics[width=0.49\textwidth]{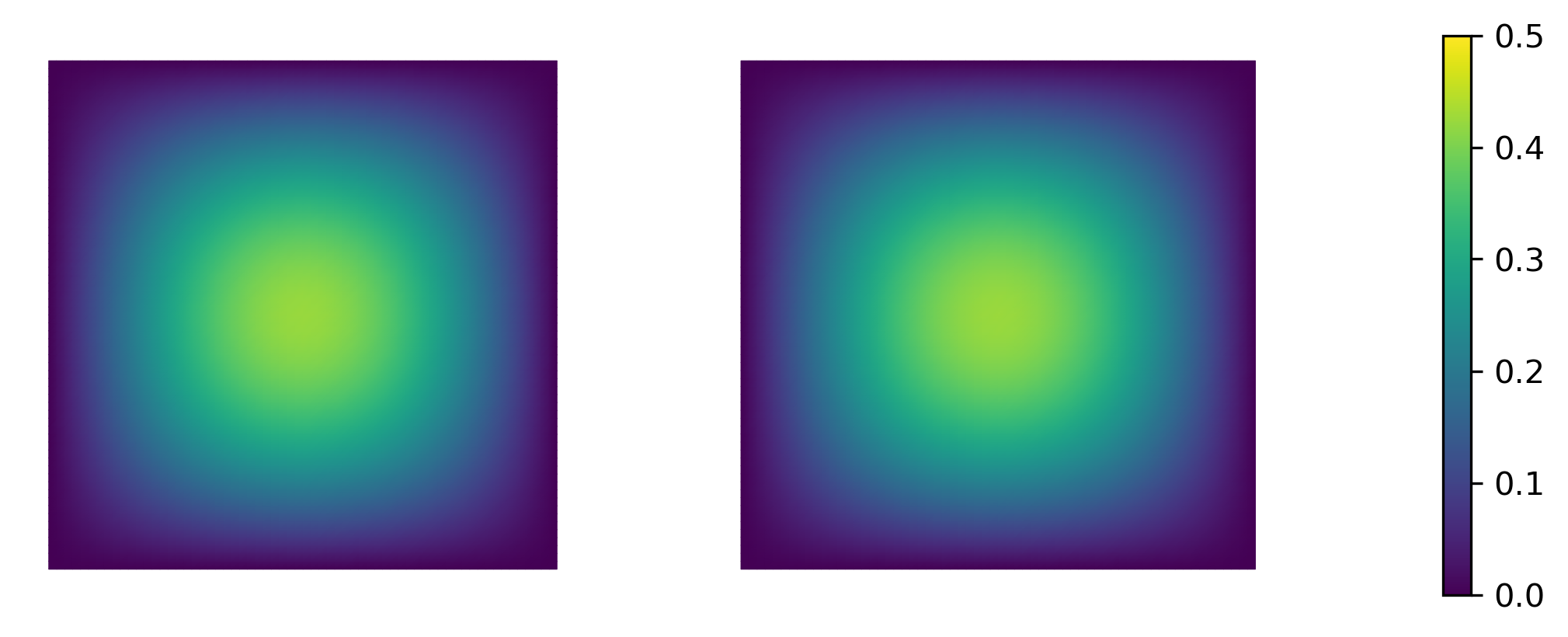}
\includegraphics[width=0.49\textwidth]{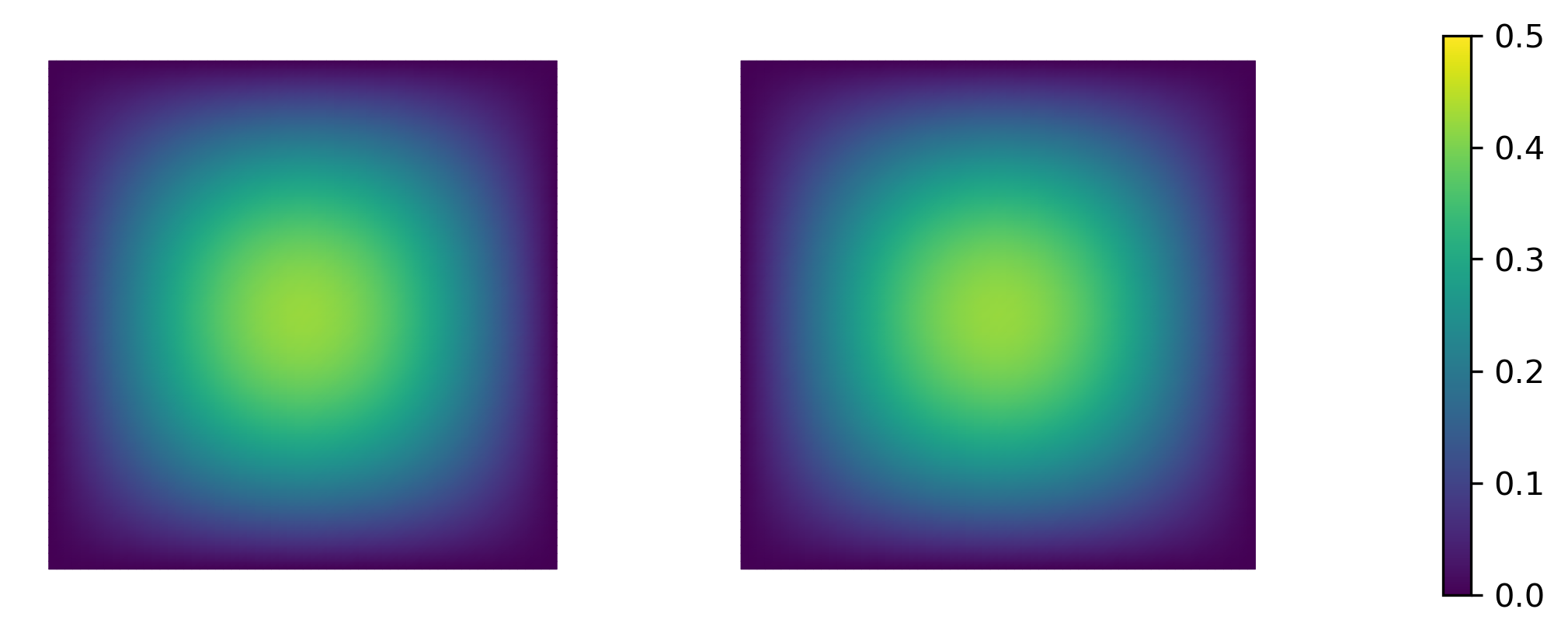}\\
\caption{Timesteps $0; 0.8; 1.6$ (left panel; top to bottom); $2.4; 3.2; 4.0$ (right panel; top to bottom)}
\end{figure}

\subsection{Galerkin discretization of the Abel integral operator}\label{sec:Galerkin}
In this section we describe another method of discretizing the Abel integral operator with the help of a Galerkin discretization. This method has been devised and investigated in \cite{vogeli:disabel}. By doing so, the ellipticity of the variational form is maintained which can be of advantage for an error analysis, as we will also see in our discrete stability estimate. We start from the equation
\begin{equation*}
I^{1-\alpha} v = g
\end{equation*}
for the Abel integral operator $I^{1-\alpha}$ cf. \eqref{eqn:Abel} or equivalently
\begin{equation}
\label{abelinteq}
\int_0^T (I^{1-\alpha} v)(s) \varphi(s) \dd{s} = \int_0^T g(s) \varphi(s) \dd{s}\mbox{ for all }\varphi \in H^{-(1-\alpha)/2}(0,T)
\end{equation}
which we abbreviate as
\begin{equation*}
a_{\alpha}(v, \varphi) = (g, \varphi) \mbox{ for all }\varphi \in H^{-(1-\alpha)/2}(0,T).
\end{equation*}
To discretize this equation we use a finite dimensional function space $V_J \subseteq H^{-(1-\alpha)/2}(0,T)$ spanned by the basis $(\phi_1,\ldots,\phi_J)$ and replace the above variational equation by
\begin{equation}
\label{abelinteq2}
a_{\alpha}(v_J, \varphi_J) = (g, \varphi_J), \quad \mbox{ for all } \varphi \in V_J
\end{equation}
where $v_J \in V_J$. We chose $V_J$ to consist of piecewise constant functions on an equidistant time grid
$t_i = i \Delta t$, $\hat{t}_{i\pm} = (i\pm\tfrac12) \Delta t$, $\hat{t}_{0-} = 0$, $\hat{t}_{J+} = T$, $\Delta t=\frac{T}{J}$, i.e., 
\begin{equation*}
v_J(t) = \sum_{i=0}^{J-1} v_i^J \phi_i(t)\,, \quad v_i^J\approx v(t_i)\,,\quad 
\phi_i(t) = \left\{ \begin{array}{cc} 1 & \mbox{ for }t \in [\hat{t}_{i-}, \hat{t}_{i+}] \\ 0 & \mbox{ else} \end{array} \right..
\end{equation*}
By choosing the test functions $\varphi_J:= \phi_j$ we can rewrite equation \eqref{abelinteq2} as
\begin{equation}\label{abelphiphij}
\frac{1}{\Gamma (1-\alpha)} \sum_{i=0}^J v_i^J \int_0^T \int_0^t (t-s)^{-\alpha} \phi_i(s) \dd{s} \, \phi_j(t) \dd{t} \overset{!}{=} \int_0^T g(t) \phi_j(t) \dd{t}\,, \quad j\in\{0,\ldots,J-1\}\,,
\end{equation}
where the left hand side can be computed as
\[ \begin{aligned}
&\frac{1}{\Gamma(1-\alpha)} \sum_{i=0}^J v_i^J \int_{\hat{t}_{j-}}^{\hat{t}_{j+}} \int_{\hat{t}_{i-}}^{\min \{ t, \hat{t}_{i+} \}} (t-s)^{-\alpha} \dd{s} \dd{t}  \\
&= \frac{1}{\Gamma (2-\alpha)} \sum_{i=0}^J v_i^J \int_{\hat{t}_{j-}}^{\hat{t}_{j+}} [(t-\hat{t}_{i-})^{1-\alpha} - (t-\min \{ t, \hat{t}_{i+} \})^{1-\alpha}] \cdot\mathbf{1}_{[\hat{t}_{i-},\infty)}(t) \dd{t} \\
&= \frac{1}{\Gamma (3- \alpha)}
\left\{ \begin{array}{ll} 
\sum_{i=0}^J v_i^J [(\hat{t}_{j+} - \hat{t}_{i-})^{2-\alpha} - ( \hat{t}_{j+} - \hat{t}_{i+})^{2-\alpha} & \\
\qquad\qquad- (\hat{t}_{j-} - \hat{t}_{i-})^{2-\alpha}+ (\hat{t}_{j-} - \hat{t}_{i+})^{2-\alpha}] & \mbox{ if } j \geq i+1 \\
\sum_{i=0}^J v_i^J (\hat{t}_{i+} - \hat{t}_{i-})^{2-\alpha} & \mbox{ if } j=i \\
0 & \mbox{ else.} 
\end{array} \right.
\end{aligned} \]
Summarizing we have
\begin{equation*}
\frac{1}{\Gamma (3- \alpha)} \sum_{i=0}^j a_{i,j}^{\alpha} v_i^J = \int_{\hat{t}_{j-}}^{\hat{t}_{j+}} g(t) \dd{t} =: \Delta t \hat{g}_j^J
\end{equation*}
with coefficients
\begin{equation*}
a_{i,j}^{\alpha} = \left\{ \begin{array}{ll} (\hat{t}_{j+} - \hat{t}_{i-})^{2-\alpha} - ( \hat{t}_{j+} - \hat{t}_{i+})^{2-\alpha} - (\hat{t}_{j-} - \hat{t}_{i-})^{2-\alpha} + (\hat{t}_{j-} - \hat{t}_{i+})^{2-\alpha} & \mbox{ if } i \leq j-1 \\
(\hat{t}_{i+}- \hat{t}_{i-})^{2-\alpha} &\mbox{ if } i=j. \end{array} \right.
\end{equation*}
That is, setting $j=n$, $l=n-i$,  $\hat{b}_l^n = (\Delta t)^{-1} a_{n-l,n}^{\alpha}$ for $l\in\{0,\ldots,n\}$, we obtain the discretization 
\begin{equation}\label{coeffGalerkin}
\begin{aligned}
&(I^{1-\alpha}v)(t_n)\approx\hat{g}_n^J = \frac{1}{\Delta t \,\Gamma (3- \alpha)} \sum_{i=0}^n a_{i,n}^{\alpha} v_i^J = \sum_{l=0}^{n} \hat{b}_l^n v(t_{n-l}) \\
&\mbox{with coefficients }\ 
\hat{b}_l^n = \frac{(\Delta t)^{1-\alpha}}{\Gamma ( 3 - \alpha)}\left\{ \begin{array}{ll} (l+1)^{2-\alpha} - 2l^{2-\alpha} +(l-1)^{2-\alpha} & \mbox{ for } l \geq 1 \\ 1 & \mbox{ for } l=0. \end{array} \right.
\end{aligned}
\end{equation}

\subsubsection{Numerical experiment with Galerkin scheme}
The numerical experiment from above is also carried out with the Galerkin type discretization of the Abel integral operator. The setting is the same, however it leads to even better results compared with the L1-approximation.

\begin{figure}[ht]
\includegraphics[width=0.49\textwidth]{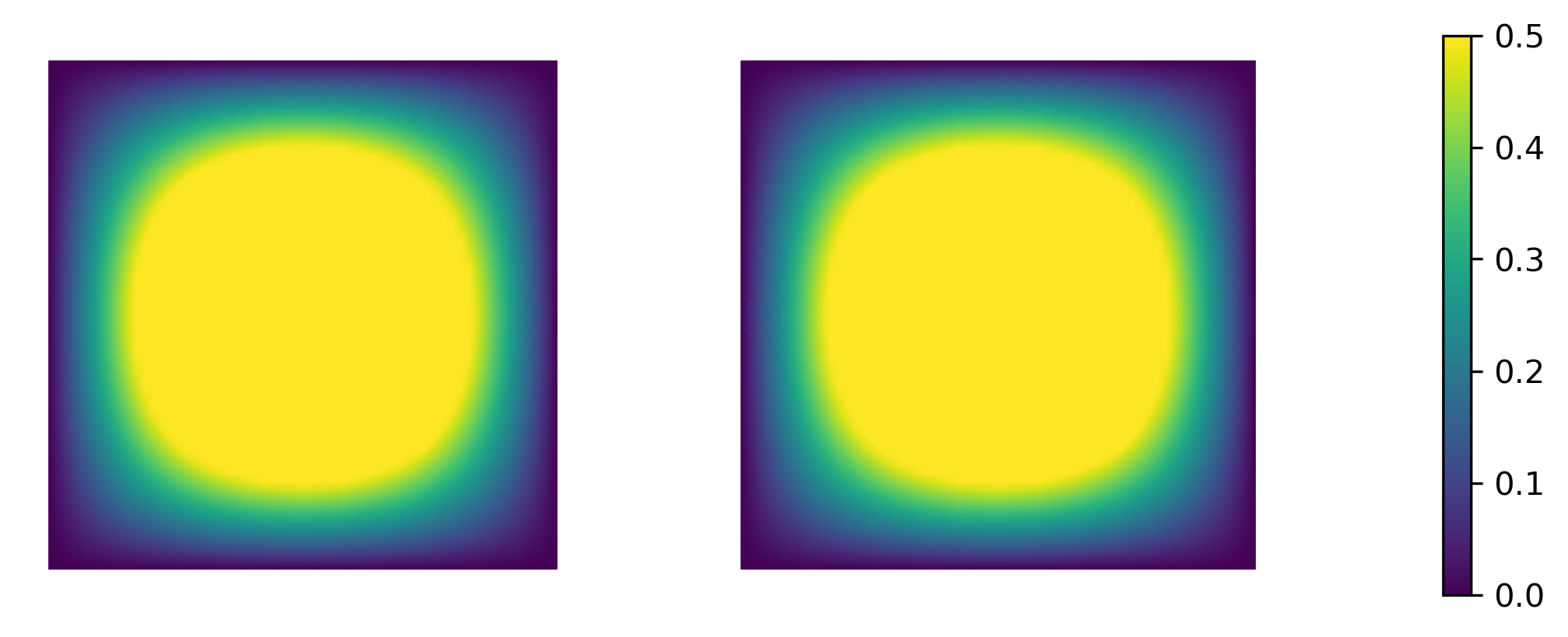}
\includegraphics[width=0.49\textwidth]{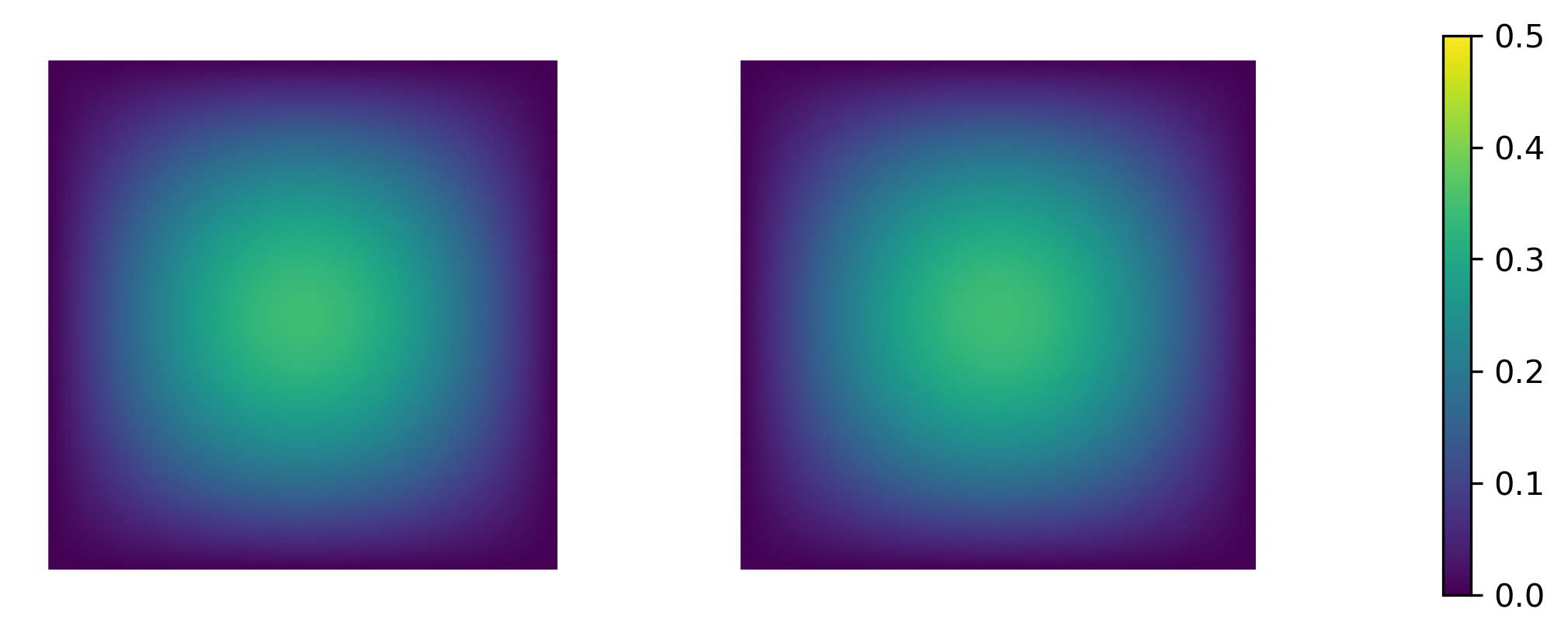}\\
\includegraphics[width=0.49\textwidth]{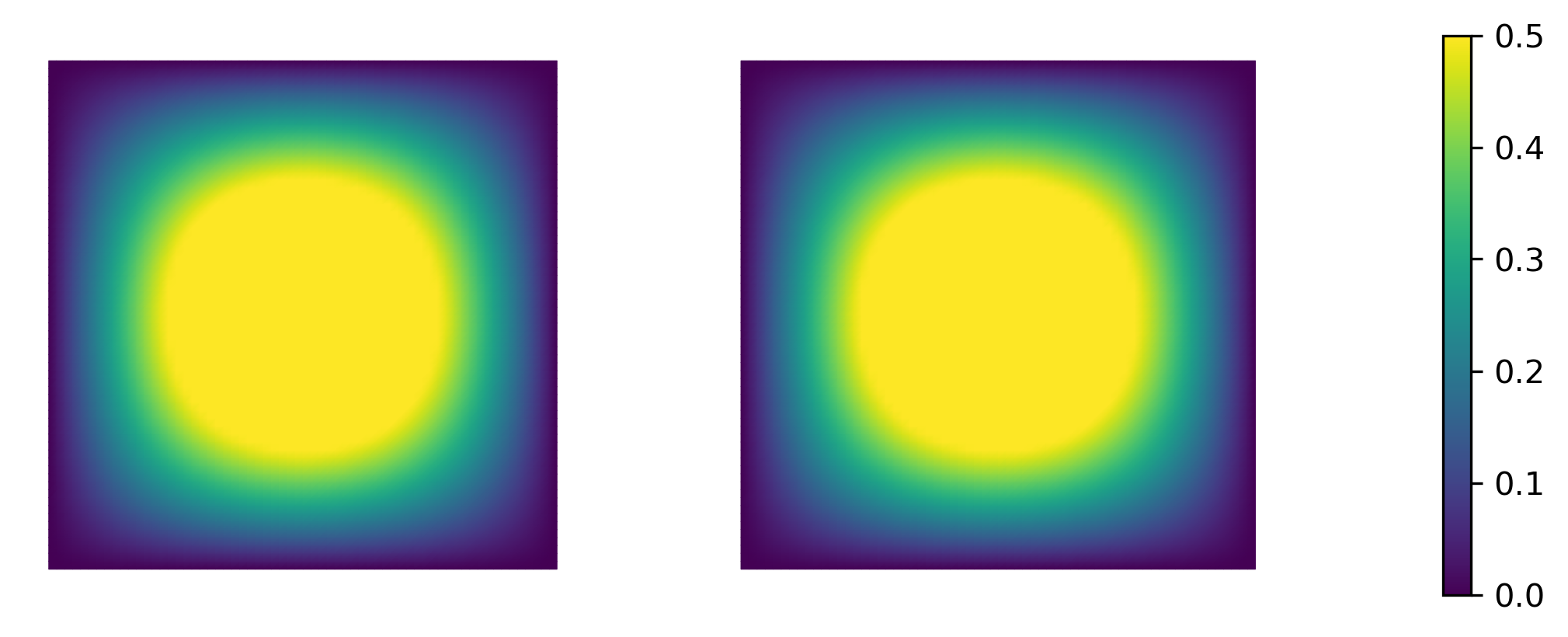}
\includegraphics[width=0.49\textwidth]{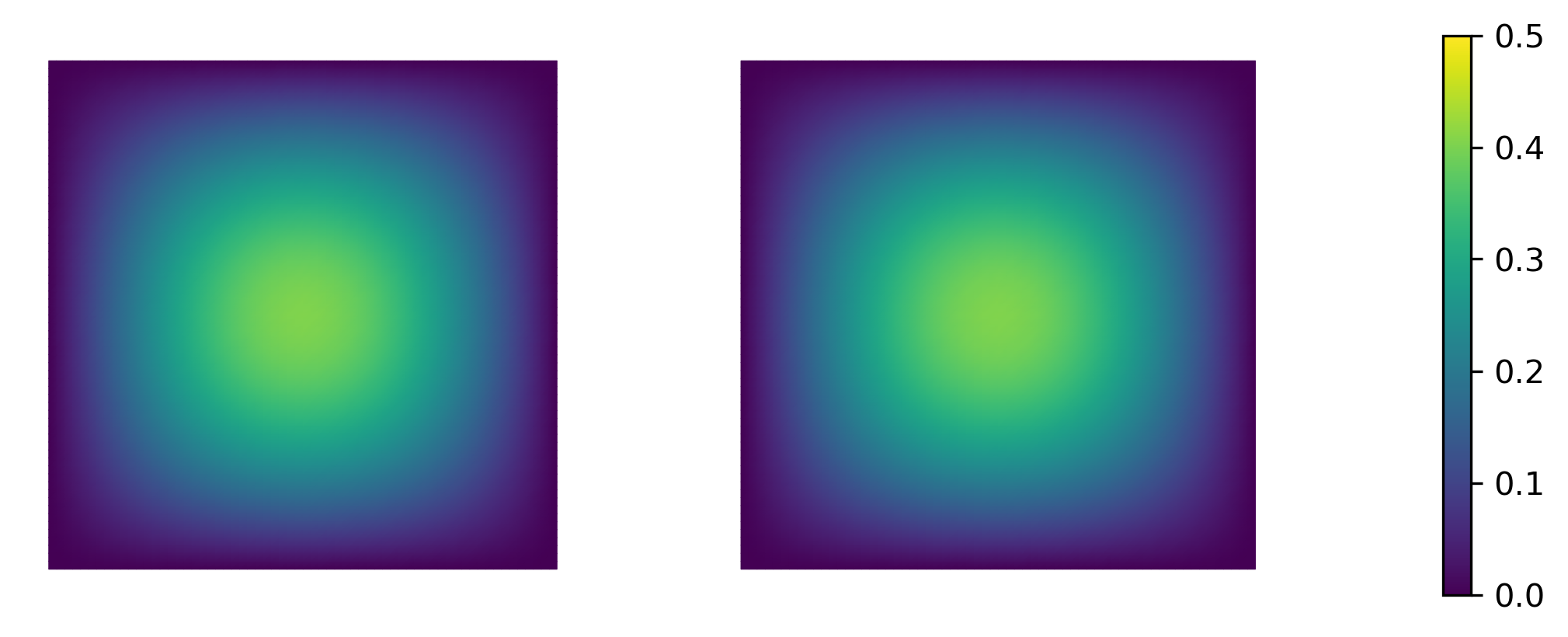}\\
\includegraphics[width=0.49\textwidth]{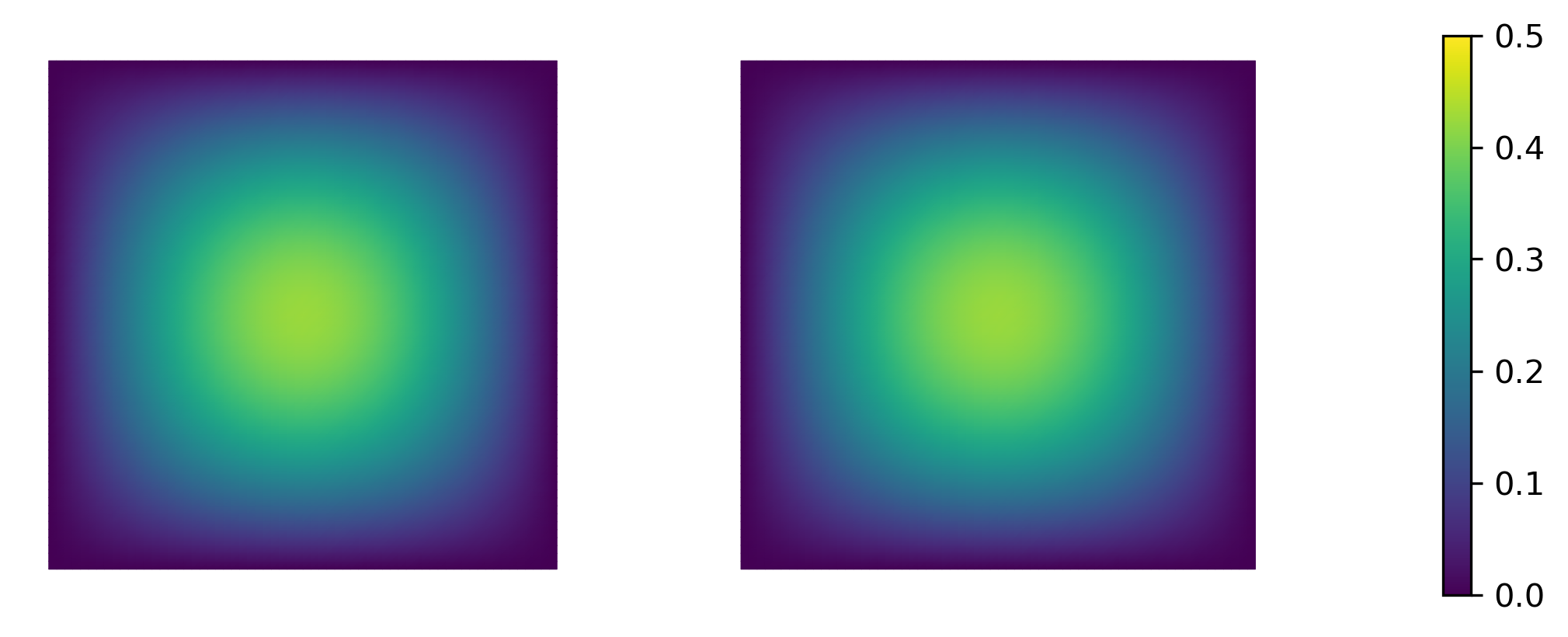}
\includegraphics[width=0.49\textwidth]{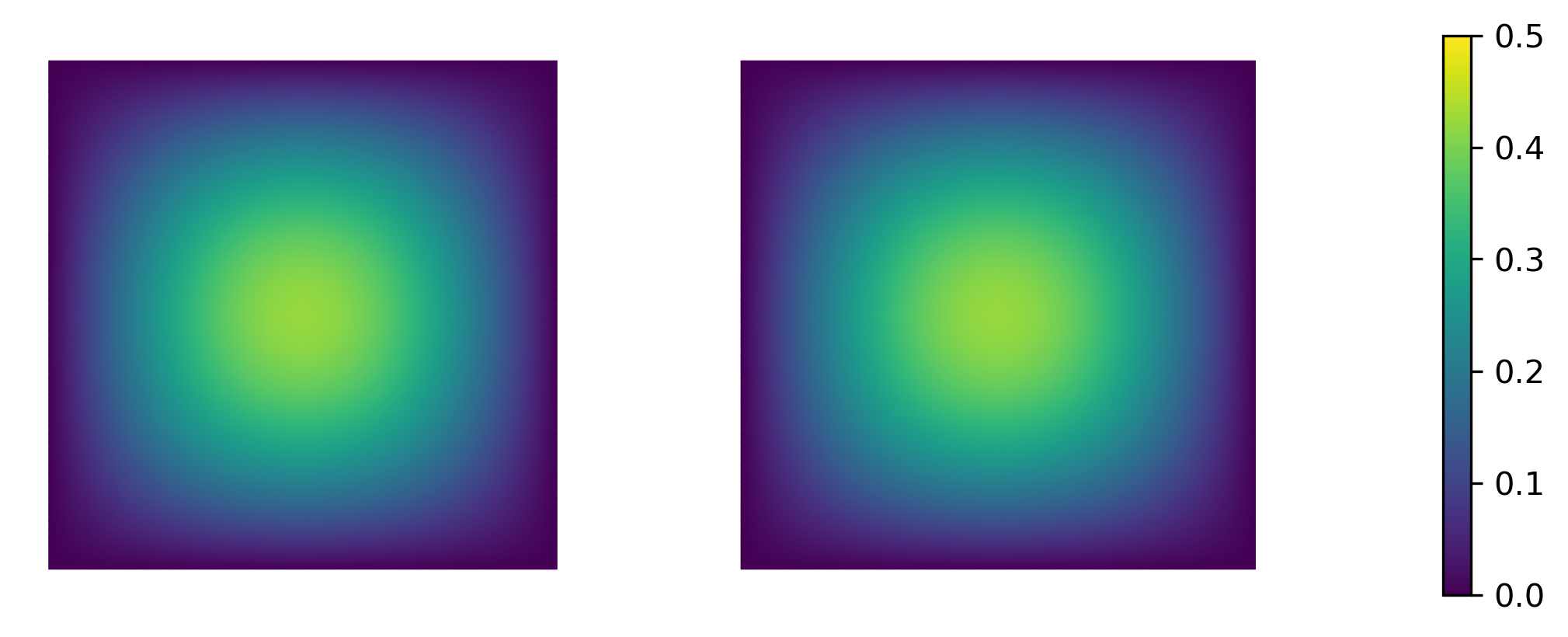}\\
\caption{
Timesteps 
$0; 0.8; 1.6$ (left panel; top to bottom); $2.4; 3.2; 4.0$ (right panel; top to bottom)
}
\end{figure}

\subsection{Error analysis}
We now carry out an error analysis of the Newmark scheme with discretization of the fractional derivative by one of the above schemes. To this end, a three stage formulation is derived, which is later used to obtain a stability estimate and convergence rates.

\subsubsection{Continuous stability}
Consider the time continuous equation \eqref{eqn:semidiscr}, rewritten by means of the Abel integral operator $I^{1-\alpha}$ as 
\begin{equation}
\label{modelcont}
M \ddot{u} + C I^{1-\alpha} \dot{u} + Ku = f,
\end{equation}
whose stability has already been shown in \cite{kaltenbacherrundell:fracpat} in the special case of a Galerkin approximation with eigenfunctions as part of the well-posedness proof of the PDE. We here recap this part of the proof with general positive semindefinite matrices $M$, $C$, and $K$, in order to later on carry it over to its time discrete version. To this end, we multiply \eqref{modelcont} with $\dot{u}$ and integrate with respect to time to obtain
\begin{equation}
\label{enestL2}
\begin{aligned}
\frac{1}{2} &| M^{1/2} \dot{u}(t)|^2 + \int_0^t \langle (I^{1-\alpha} C^{1/2} \dot{u})(s), C^{1/2} \dot{u}(s) \rangle \dd{s} + \frac{1}{2} |K^{1/2} u(t)|^2 \\
&\leq \frac{1}{2}|M^{1/2} \dot{u}(0)|^2 + \frac{1}{2} |K^{1/2} u(0)|^2 + \frac{1}{2} \int_0^t |M^{-1/2}f|^2 \dd{s} + \frac{1}{2} \int_0^t |M^{1/2} \dot{u} |^2 \dd{s}
\end{aligned}
\end{equation}
or alternatively
\begin{equation}
\label{enestL1}
\begin{aligned}
\frac{1}{2} &| M^{1/2} \dot{u}(t)|^2 + \int_0^t \langle (I^{1-\alpha} C^{1/2} \dot{u})(s), C^{1/2} \dot{u}(s) \rangle \dd{s} + \frac{1}{2} |K^{1/2} u(t)|^2 \\
&\leq \frac{1}{2}|M^{1/2} \dot{u}(0)|^2 + \frac{1}{2} |K^{1/2} u(0)|^2 + 
\Bigl(\int_0^t |M^{-1/2}f| \dd{s}\Bigr)^2 + \frac{1}{4} \sup_{s\in[0,t]} |M^{1/2} \dot{u}(s) |^2 
\end{aligned}
\end{equation}

From \cite{vogeli:disabel} we use the estimate
\begin{equation}
\label{corecivityI}
\int_0^t \langle (I^{1-\alpha} C^{1/2} \dot{u})(s), C^{1/2} \dot{u}(s) \rangle \dd{s} \geq \cos ( \tfrac{\pi(1-\alpha)}{2} ) \| C^{1/2} \dot{u} \|_{H^{-(1-\alpha)/2}(0,t)}^2 \geq 0
\end{equation} 
and with Gronwall's inequality we obtain from \eqref{enestL2}
\begin{equation}
\label{stabestL1}
\begin{aligned}
&|M^{1/2} \dot{u}(t) |^2 + \| C^{1/2} \dot{u} \|_{H^{-(1-\alpha)/2}(0,t)}^2 + |K^{1/2} u(t)|^2 \\
&\leq c(T) \left( |M^{1/2} \dot{u}(0) |^2 + |K^{1/2} u(0)|^2 + \int_0^t |M^{-1/2} f |^2 \dd{s}\right),
\end{aligned}
\end{equation}
and alternatively from \eqref{enestL1}
\begin{equation}
\label{stabestL2}
\begin{aligned}
&\sup_{s\in[0,t]}|M^{1/2} \dot{u}(s) |^2+\| C^{1/2} \dot{u} \|_{H^{-(1-\alpha)/2}(0,t)}^2+|K^{1/2} u(t)|^2\\
&\leq c \left( |M^{1/2} \dot{u}(0) |^2 + |K^{1/2} u(0)|^2 + \Bigl(\int_0^t |M^{-1/2} f | \dd{s}\Bigr)^2\right),
\end{aligned}
\end{equation}
for some (possibly large) constants $c,c(T)>0$. 
\subsubsection{Three stage formulation}
To obtain a three stage formulation for the discrete problem, analogously to, e.g., \cite{EhrlicherBonaventuraBursi2002}, we first of all rewrite the discrete model \eqref{Newmark} with $u$ and its derivatives written as displacement $d=u$, velocity $v=\dot{u}$ and acceleration $a=\ddot{u}$ (note that for convenience we are using mechanics notation in this section, while in our PAT application, these quantities clearly have different physical meanings)
\begin{equation}
\label{modeldis}
M a_{n+1} + C \sum_{j=0}^{n+1} b_j^{n+1} v_{n+1-j} + K d_{n+1} = f_{n+1}\,.
\end{equation}
For simplicity of exposition and since it is the most favorable case in terms of consistency, we focus on the Newmark parameters $\gamma = \frac{1}{2}$ and $\beta = \frac{1}{4}$. We aim for a formulation which only depends on $v$ and for this purpose use the identities
\begin{equation}
\label{vnpdnp}
\begin{aligned}
v_{n+1} &= v_n + \frac{\Delta t}{2}(a_n + a_{n+1}) \\
d_{n+1} &= d_n + \Delta t v_n + \frac{\Delta t^2}{2} (a_n + a_{n+1}) = d_n + \frac{\Delta t}{2}(v_n + v_{n+1})
\end{aligned}
\end{equation}
Now we sum $\frac{1}{4} (\ref{modeldis})_{n+1} + \frac{1}{2} (\ref{modeldis})_{n} + \frac{1}{4} (\ref{modeldis})_{n-1}$ to obtain
\begin{equation}\label{modeldis_lincomb}
\begin{aligned}
&M \left[ \frac{v_{n+1}-v_{n-1}}{2 \Delta t} \right] + C \left[ \underbrace{ \frac{1}{4} \sum_{j=0}^{n+1} b_j^{n+1} v_{n+1-j} + \frac{1}{2} \sum_{i=0}^n b_i^n v_{n-i} + \frac{1}{4} \sum_{k=0}^{n-1} b_k^{n-1} v_{n-1-k}}_{=: \tilde{v}_n^{\alpha}} \right] \\
&+ K \left[ \underbrace{ \frac{d_{n+1} + 2 d_n + d_{n-1}}{4} }_{=:\tilde{d}_n}\right] = \underbrace{ \frac{f_{n+1} +2f_n +f_{n-1}}{4} }_{=:\tilde{f}_n}
\end{aligned}
\end{equation}
As it holds $b_j^{n+1} = b_j^{n} = b_j^{n-1}$ for $0 \leq j \leq n-1$ we obtain for $\tilde{v}_n^{\alpha}$
\begin{equation*}
\begin{aligned}
&\tilde{v}_n^{\alpha} = \sum_{j=0}^n b_j^n \tilde{v}_{n-j} + \frac{b_{n+1}^{n+1}v_0+b_n^n v_1}{4}\\
&\mbox{ with }\tilde{v}_{i}=\frac{v_{i+1} +2v_i + v_{i-1}}{4}, i\in\{1,\ldots,n\}, \quad
\tilde{v}_0 = \frac{v_1+v_0}{2}\,.
\end{aligned}
\end{equation*}
For the displacement $\tilde{d}_n$ we obtain a new expression by considering the difference $\tilde{d}_n - \tilde{d}_{n-1}$
\begin{equation*}
\begin{aligned}
\tilde{d}_n - \tilde{d}_{n-1} &= \frac{1}{4}( d_{n+1} - d_n + 2(d_{n} - d_{n-1}) + d_{n-1} -d_{n-2})\\ 
&= \Delta t \frac{v_{n+1} + 3v_n + 3v_{n-1} + v_{n-2}}{8} =: \Delta t \tilde{\tilde{v}}_n
\end{aligned}
\end{equation*}
where we used that $d_n - d_{n-1} = \frac{\Delta t}{2} (v_n + v_{n-1})$. By rearrangement we obtain for $\tilde{d}_n$
\begin{equation*}
\tilde{d}_n = \tilde{d}_1 + \Delta t \sum_{j=2}^n \tilde{\tilde{v}}_j = d_0 + \Delta t \sum_{j=1}^n \tilde{\tilde{v}}_j 
\end{equation*}
with $\tilde{\tilde{v}}_1 := \frac{v_2 + 4v_1 + 3v_0}{8}$. 
Using the fact that $\tilde{\tilde{v}}_j = \frac{\tilde{v}_j+\tilde{v}_{j-1}}{2}$ for $j\in\{0,\ldots,n\}$, we arrive at the following three stage formulation for the discretized model
\begin{equation}
\label{eqdisc}
M \left[ \frac{v_{n+1} - v_{n-1}}{2 \Delta t} \right] + C \left[ \sum_{j=0}^n b_j^{n+1} \tilde{v}_{n-j} + 
\frac{b_{n+1}^{n+1}v_0+b_n^n v_1}{4}
\right] + K \left[ d_0 + \Delta t \sum_{j=1}^n \frac{\tilde{v}_j + \tilde{v}_{j-1}}{2} \right] = \tilde{f}_n
\end{equation}

\subsubsection{Stability for the discrete model}
To derive a discrete stability equation we will proceed analogously to the time continuous case \eqref{modelcont}, that is, multiply its discretized counterpart \eqref{eqdisc} by $\Delta t \tilde{v}_n$ and sum up over $n$.

To obtain estimates from below for the resulting left hand side, we need the following two Lemmas.

\begin{lemma}
\label{lemma1}
For any vector $w_0,\ldots,w_{N+1}$
\begin{equation}
\sum_{n=1}^N \left(w_{n+1} - w_{n-1}, \frac{w_{n+1} + 2 w_n + w_{n-1}}{4}\right) = \left| \frac{w_{N+1} + w_N}{2} \right|^2 - \left|\frac{w_1 + w_0}{2} \right|^2
\end{equation}
\end{lemma}

\textit{Proof.} see the appendix.
\begin{flushright}
$\diamondsuit$
\end{flushright}

\begin{lemma}
\label{lemma2}
For any vector $\tilde{w}_0,\ldots,\tilde{w}_{N}$
\begin{equation}
\sum_{n=0}^N \left(\sum_{j=1}^n \frac{\tilde{w}_j + \tilde{w}_{j-1}}{2}, \tilde{w}_n \right) \geq \frac{1}{4} \left[ \left| \sum_{n=0}^N \tilde{w}_n \right|^2 - \left| \tilde{w}_0 \right|^2 \right]
\end{equation}
\end{lemma}

\textit{Proof.} see the Appendix.
\begin{flushright}
$\diamondsuit$
\end{flushright}

Now, we can consider the inner product of (\ref{eqdisc}) with $\Delta t\, \tilde{v}_n$, and sum up over $n$ to obtain the following inequalities from Lemmas \ref{lemma1} and \ref{lemma2} 
\[ \begin{aligned}
&\Delta t \sum_{n=1}^N\left( M \frac{v_{n+1}-v_{n-1}}{2 \Delta t}, \frac{v_{n+1} + 2 v_n + v_{n-1}}{4} \right) \geq \frac{1}{2} \left[ |M^{1/2} \frac{v_{N+1}+v_N}{2} |^2 - |M^{1/2} \frac{v_1 +v_0}{2} |^2 \right] \\
&\Delta t\sum_{n=0}^N \left( K (d_0 + \Delta t \sum_{j=1}^n \frac{\tilde{v}_j+\tilde{v}_{j-1}}{2} ), \tilde{v}_n \right)  \\
&\hspace*{4cm}\geq ( K^{1/2} d_0, K^{1/2} \Delta t \sum_{n=0}^N \tilde{v}_n ) + \frac{1}{4} |K^{1/2} \Delta t \sum_{n=0}^N \tilde{v}_n |^2 - \frac{1}{4} | K^{1/2} \Delta t \tilde{v}_0|^2 \\
\\
&\hspace*{4cm}\geq -2|K^{1/2} d_0|^2 + \frac{1}{8} |K^{1/2} \Delta t \sum_{n=0}^N \tilde{v}_n |^2 - \frac{1}{4} | K^{1/2} \Delta t \tilde{v}_0|^2 \\
&\Delta t\sum_{n=0}^N \left( C \left( \sum_{j=0}^n b_j^n \tilde{v}_{n-j} + r_0^n \right), \tilde{v}_n \right)  \\
&\hspace*{4cm}= \int_0^{(N+\frac12)\Delta t} C^{1/2} [(I^{1-\alpha} \tilde{v}_J)(t)+r^0_J(t)], C^{1/2} v_J(t)) \dd{t} \\
&\hspace*{4cm}\geq \frac{\cos ( \tfrac{\pi(1-\alpha)}{2} )}{2} \| C^{1/2} \tilde{v}_J \|_{H^{-(1-\alpha)/2}(0,t)}^2 - \frac{1}{2\cos ( \tfrac{\pi(1-\alpha)}{2} )} \| C^{1/2} r^0_J \|_{H^{(1-\alpha)/2}(0,t)}^2,
\end{aligned} \]
where the latter holds in case of the Galerkin discretization and we have used the abbreviations
\begin{equation}
\label{defr0J}
r_0^n=\frac{b_{n+1}^{n+1}v_0+b_n^n v_1}{4}\,, \quad \tilde{v}_J(t)=\sum_{n=0}^{J-1} \tilde{v}_n\phi_n(t)
\,, \quad r^0_J(t)=\sum_{n=0}^{J-1} r^0_n\phi_n(t)\,.
\end{equation}

Finally, to estimate the right hand side in terms of the mass term (as we have done in the continuous setting), we abbreviate $\hat{v}_n=\frac{v_{n+1}+v_n}{2}$ and use the identity $\tilde{v}_n=\frac{\hat{v}_n+\hat{v}_{n-1}}{2}$ as well as
\[ \begin{aligned}
&\Delta t\sum_{n=1}^N \left(\tilde{f}_n, \tilde{v}_n\right) 
= \Delta t\sum_{n=1}^N \left(M^{-1/2}\tilde{f}_n, M^{1/2} \frac{\hat{v}_n+\hat{v}_{n-1}}{2}\right) \\
&= \Delta t\sum_{n=1}^{N-1} \left(M^{-1/2}\frac{\tilde{f}_{n+1}+\tilde{f}_n}{2}, M^{1/2} \hat{v}_n\right) 
+\frac{\Delta t}{2} \left(M^{-1/2}\tilde{f}_1, M^{1/2} \hat{v}_0\right)
+\frac{\Delta t}{2} \left(M^{-1/2}\tilde{f}_N, M^{1/2} \hat{v}_N\right)\\
&\leq \frac{\Delta t}{2}\left(\sum_{n=1}^N\left| M^{-1/2}\frac{\tilde{f}_{n+1}+\tilde{f}_n}{2}\right|^2 
+ \frac14 \left| M^{-1/2}\tilde{f}_1\right|^2+ \frac14 \left| M^{-1/2}\tilde{f}_N\right|^2\right)
+ \frac{\Delta t}{2} \sum_{n=0}^N\left|M^{1/2} \hat{v}_n\right|^2
\end{aligned} \]
by Cauchy-Schwarz' and Young's inequalities, or alternatively
\[ \begin{aligned}
&\Delta t\sum_{n=1}^N \left(\tilde{f}_n, \tilde{v}_n\right) \\
&\leq \Delta t^2\left(\sum_{n=1}^N\left| M^{-1/2}\frac{\tilde{f}_{n+1}+\tilde{f}_n}{2}\right| 
+ \frac12 \left| M^{-1/2}\tilde{f}_1\right|+ \frac12 \left| M^{-1/2}\tilde{f}_N\right|\right)^2
+ \frac14 \max_{n\in\{1\,\ldots,N\}}\left|M^{1/2} \hat{v}_n\right|^2
\end{aligned} \]

Combining these estimates with the discrete Gronwall inequality
\[
\begin{aligned}
&\eta^{(n)} \leq a^{(n)} +\Delta t\sum_{j=1}^n b^{(j)} \eta^{(j)} \quad\mbox{ for all }n\in\{1,\ldots, N\}
\\&\ \Rightarrow \
\eta^{(n)} \leq a^{(n)} +\Delta t\sum_{j=1}^n a^{(j)}\, b^{(j)} \, \exp\Bigl(\Delta t\sum_{i=j}^n b^{(i)}\Bigr) \quad
\mbox{ for all }n\in\{1,\ldots, N\}\,,
\end{aligned}
\] 
we obtain the following discrete stability estimate for a solution to \eqref{eqdisc}.
\begin{theorem}\label{th:stbdiscr}
For the time discretization scheme \eqref{eqdisc} with coefficients $b_l^n=\hat{b}_l^n$ according to \eqref{coeffGalerkin} there exists a constant $C(T)$ depending only on $T$ such that with $\hat{v}_n=\frac{v_{n+1}+v_n}{2}$, $n\in\{0,\ldots,J\}$, $\tilde{v}_n= \frac{v_{n+1} +2v_n + v_{n-1}}{4}$, $n\in\{1,\ldots,J\}$, $J=\frac{T}{\Delta t}$, $\tilde{v}_0 = \frac{v_1+v_0}{2}$, $\tilde{v}_n=\sum_{n=0}^{J-1} \tilde{v}_n \phi_n(t)$, $r^0_J$ as in \eqref{defr0J},
\begin{equation}\label{stability_discreteL2}
\begin{aligned}
&|M^{1/2} \hat{v}_N |^2 + \| C^{1/2} \tilde{v}_J \|_{H^{-(1-\alpha)/2}(0,t)}^2 + |K^{1/2} \Delta t \sum_{n=0}^N \tilde{v}_n |^2\\
&\leq 
c(T)\Bigl(|M^{1/2} \hat{v}_0 |^2 + | K^{1/2} d_0|^2 +  | \Delta t K^{1/2} \tilde{v}_0|^2 
+\| C^{1/2} r^0_J \|_{H^{(1-\alpha)/2}(0,t)}^2
+ \sum_{n=0}^N\left| M^{-1/2}\tilde{f}_n\right|^2\Bigr)
\end{aligned}
\end{equation}
\begin{equation}\label{stability_discreteL1}
\begin{aligned}
&\max_{n\in\{1\,\ldots,N\}} |M^{1/2} \hat{v}_n |^2 + \| C^{1/2} \tilde{v}_J \|_{H^{-(1-\alpha)/2}(0,t)}^2 
+ |K^{1/2} \Delta t \sum_{n=0}^N \tilde{v}_n |^2\\
&\leq 
c\Bigl(|M^{1/2} \hat{v}_0 |^2  + | K^{1/2} d_0|^2 +  |\Delta t K^{1/2} \tilde{v}_0|^2 
+\| C^{1/2} r^0_J \|_{H^{(1-\alpha)/2}(0,t)}^2+ 
\Bigl(\sum_{n=0}^N\left| M^{-1/2}\tilde{f}_n\right|\Bigr)^2\Bigr)
\end{aligned}
\end{equation}
for all $N\in\{1,\ldots,\frac{T}{\Delta t}\}$.
\end{theorem}
Estimates \eqref{stability_discreteL2}, \eqref{stability_discreteL1} are the time discrete counterparts of the continuous stability estimates \eqref{stabestL1}, \eqref{stabestL2}, respectively.

\subsubsection{Global error}
For the estimation of the global error consider the identity
\begin{equation*}
 \frac{1}{2\Delta t} \int_{t_{n-1}}^{t_{n+1}} \ddot{u}(t) \dd{t} = \frac{v(t_{n+1})-v(t_{n-1})}{2\Delta t}
\end{equation*}
so that after integrating, $\int_{t_{n-1}}^{t_{n+1}} \dd{t}$, the continuous system \eqref{eqn:semidiscr} can be reformulated analogously to \eqref{eqdisc} as 
\begin{equation}\label{eqn:semidiscr_1}
M \frac{v(t_{n+1}) -v(t_{n-1})}{2\Delta t} + C \frac{1}{2\Delta t} \int_{t_{n-1}}^{t_{n+1}} (I^{1-\alpha} v)(t) \dd{t} + K \frac{1}{2\Delta t} \int_{t_{n-1}}^{t_{n+1}} u(t) \dd{t} = \frac{1}{2\Delta t} \int_{t_{n-1}}^{t_{n+1}} f(t) \dd{t}.
\end{equation}
We denote the error by
\begin{equation}
\label{en}
\dot{e}_n = v(t_n)-v_n = \dot{u}(t_n)-\dot{u}_n\,,  \quad
e_n = u(t_n)-u_n
\end{equation}
and the related quantities
\begin{equation}
\label{enhat}
\hat{\dot{e}}_i=\frac{\dot{e}_{i+1}+\dot{e}_i}{2}\,, \quad
\tilde{\dot{e}}_i = \frac{\dot{e}_{i+1} + 2\dot{e}_i + \dot{e}_{i-1}}{4}=\frac{\hat{\dot{e}}_i+\hat{\dot{e}}_{i-1}}{2}\,, \quad 
\tilde{\dot{e}}_0 = \frac{\dot{e}_1+\dot{e}_0}{2}
\end{equation}
analogously to the previous subsection. 
Subtracting \eqref{eqdisc}  from \eqref{eqn:semidiscr_1} and making use of the identity of \eqref{modeldis_lincomb} and \eqref{eqdisc} with $v_i$, $d_i$ replaced by $v(t_i)$, $u(t_i)$, respectively, we obtain the following
\begin{equation}\label{eqne}
\begin{aligned}
&M \frac{\dot{e}_{n+1} - \dot{e}_{n-1}}{2 \Delta t}  + C \left(\sum_{j=0}^n b_j^n \tilde{\dot{e}}_{n-j}+ \frac{b_{n+1}^{n+1}\dot{e}_0+b_n^n \dot{e}_1}{4}\right) 
+ K \left( \Delta t \sum_{j=1}^n \frac{\tilde{\dot{e}}_j + \tilde{\dot{e}}_{j-1}}{2} \right) \\
&= C \left[ \frac{1}{4} \sum_{j=0}^{n+1} b_j^{n+1} v(t_{n+1-j}) + \frac{1}{2} \sum_{i=0}^n b_i^n v(t_{n-i}) + \frac{1}{4} \sum_{k=0}^{n-1} b_k^{n-1} v(t_{n-1-k}) 
- \frac{1}{2\Delta t} \int_{t_{n-1}}^{t_{n+1}} (I^{1-\alpha} v)(t) \dd{t}\right] \\
&\quad+ K \left[ \frac{u(t_{n+1}) + 2 u(t_n) + u(t_{n-1})}{4} 
- \frac{1}{2\Delta t} \int_{t_{n-1}}^{t_{n+1}} u(t) \dd{t}\right]\\
&\quad +\frac{1}{2 \Delta t} \int_{t_{n-1}}^{t_{n+1}} f(t) \dd{t} - \frac{f_{n+1}+2f_n+f_{n-1}}{4},
\end{aligned}
\end{equation}
that is, \eqref{eqdisc} with $d_0=0$, $v$ replaced by $e$, and a different right hand side. Hence we are in the position to apply Theorem \ref{th:stbdiscr}, after estimating the discrete $L^1(0,T)$ norm of the right hand side in \eqref{eqne}.

This can be done as follows.

The term multiplied with $K$ by Taylor's Theorem can be written as
\[
\begin{aligned}
&K\Bigl(\frac{u(t_{n+1}) + 2 u(t_n) + u(t_{n-1})}{4} 
- \frac{1}{2\Delta t} \int_{t_{n-1}}^{t_{n+1}} u(t) \dd{t}\Bigr)\\
&= K\Bigl(u(t_n)+\frac14 \ddot{u}(\xi)\Delta t^2 
-\frac{1}{2\Delta t}\int_{t_{n-1}}^{t_{n+1}} 
u(t_n)+ \dot{u}(t_n)(t-t_n)+ \frac12 \ddot{u}(\zeta)(t-t_n)^2 \dd{t}\Bigr)\\
\end{aligned}
\]
so that
\begin{equation}\label{estKterm}
\left|M^{-1/2}K\Bigl(\frac{u(t_{n+1}) + 2 u(t_n) + u(t_{n-1})}{4} 
- \frac{1}{2\Delta t} \int_{t_{n-1}}^{t_{n+1}} u(t) \dd{t}\Bigr)\right|
= \|M^{-1/2}K\ddot{u}\|_{L^\infty(0,T)}O(\Delta t^2)\,. 
\end{equation}

For the term multiplied with $C$ we mention the error estimate in \cite[Theorem 21]{vogeli:disabel} that is based on $H^{-(1-\alpha)/2}(0,T)$ ellipticity of the Abel integral operator and C\'{e}a's Lemma. However, since we need an estimate in the discrete $L^1(0,T)$ norm here, we cannot rely on this result but have to estimate by means of Taylor expansion here too. 
To this end, we make use of the fact that by the derivation in Section \ref{sec:Galerkin} \eqref{abelphiphij}--\eqref{coeffGalerkin}, as well as Galerkin orthogonality and setting 
$v_{cont,J}(t)=\sum_{n=0}^{J-1} v(t_n)\phi_n(t)$, $V^\alpha(t)= (I^{1-\alpha} v)(t)=\partial_t^\alpha u(t)$ this term can be written as 
\[
\begin{aligned}
&\frac{1}{\Delta t} \int_0^T 
(I^{1-\alpha} v_{cont,J})(t)\Bigl(\frac14\phi_{n+1}(t)+\frac12\phi_n(t)+\frac14\phi_{n-1}(t)\Bigr)\dd{t}
-\frac12\int_{t_{n-1}}^{t_{n+1}} (I^{1-\alpha} v)(t)\dd{t}\\
&=\frac{1}{\Delta t} \int_0^T 
\Bigl(\frac14\phi_{n+1}(t)+\frac12\phi_n(t)+\frac14\phi_{n-1}(t)
-\frac12 1_{[t_{n-1},t_{n+1}]} (I^{1-\alpha} v)(t)\Bigr)\dd{t}\\
&=\frac{1}{4\Delta t}  
\Bigl(
\int_{t_{n+1}}^{\hat{t}_{n+1+}} V^\alpha(t)\dd{t} -\int_{\hat{t}_{n+}}^{t_{n+1}} V^\alpha(t)\dd{t}
+\int_{\hat{t}_{n-1-}}^{t_{n-1}} V^\alpha(t)\dd{t}-\int_{\hat{t}_{n-}}^{t_{n-1}} V^\alpha(t)\dd{t} \Bigr)
\end{aligned}
\]
and therefore 
\begin{equation}\label{estCterm}
\begin{aligned}
&\left|M^{-1/2}C \left[ \frac{1}{4} \sum_{j=0}^{n+1} b_j^{n+1} v(t_{n+1-j}) + \frac{1}{2} \sum_{i=0}^n b_i^n v(t_{n-i}) + \frac{1}{4} \sum_{k=0}^{n-1} b_k^{n-1} v(t_{n-1-k}) 
- \frac{1}{2\Delta t} \int_{t_{n-1}}^{t_{n+1}} (I^{1-\alpha} v)(t) \dd{t}\right]\right|\\
&=\frac{1}{96}\|M^{-1/2}C\ddot{V}^\alpha\|_{L^\infty(0,T)} \Delta t^2
\end{aligned}
\end{equation}
where we have used the Taylor expansion $V^\alpha(t)=V^\alpha(a)+\dot{V}^\alpha(a)(t-a)+\frac12\ddot{V}^\alpha(\xi_t)(t-a)^2$ for $a=t_{n\pm1}$ and some $\xi_t\in[a,t]$ or $[t,a]$, respectively, as well as cancellation of the $\dot{V}^\alpha$ term due to symmetry.

To remove the error due to $f$ in the right hand side of \eqref{eqne}, we might aim at choosing $f_n$ such that
\begin{equation}
\frac{1}{2 \Delta t} \int_{t_{n-1}}^{t_{n+1}} f(t) \dd{t} = \frac{f_{n+1} + 2 f_n + f_{n-1}}{4}
\end{equation}
which via $\frac{f_{n+1}+f_n}{2} = \frac{1}{\Delta t} \int_{t_n}^{t_{n+1}} f(t) \dd{t}$ can be achieved by setting
\begin{equation}
\quad f_{n+1} = \frac{1}{\Delta t} \sum_{j=0}^n (-1)^{n-1} \int_{t_j}^{t_{j+1}} f(t) \dd{t}.
\end{equation}
Alternatively, when simply setting $f_{n+1}=f(t_{n+1})$, this error can be estimated as above by 
\begin{equation}\label{estfterm}
\left|M^{-1/2}\Big[\frac{1}{2 \Delta t} \int_{t_{n-1}}^{t_{n+1}} f(t) \dd{t} - \frac{f_{n+1}+2f_n+f_{n-1}}{4}\Big]\right|
=\|M^{-1/2}\ddot{f}\|_{L^\infty(0,T)} O(\Delta t^2)
\end{equation}

It only remains to estimate the terms pertaining to initial conditions in \eqref{stability_discreteL1}.
From \eqref{solve_mass}, \eqref{corr_mass}, \eqref{modelcont}  as well as $\ddot{u}_0=\ddot{u}(t_0)$, $\dot{u}_0=\dot{u}(t_0)$, $u_0=u(t_0)$, we obtain
\[
\tilde{\dot{e}}_0=\hat{\dot{e}}_0=\frac12 \dot{e}_1 = \frac12 (v(t_1)-v_1) = \frac{\gamma\Delta t}{2}(M^*)^{-1}\Bigl(
C\bigl(b_0^1v(t_1)+b_1^1v(t_0)-(I^{1-\alpha}v)(t_1)\bigr)+f(t_1)-f_1\Bigr)\,,
\]
so that for we can estimate
\[
|A\dot{e}_1|=\|A(M^*)^{-1}C\frac{d^k}{dt^k}V^\alpha\|_{L^\infty(0,T)} \, O(\Delta t^{k+1})\,, 
\quad A\in \{M^{1/2},K^{1/2},CM^{1/2}\}\,, \quad k\in\{0,1,2\}\,.
\]

With $r_0^n=\frac{b_n^n \dot{e}_1}{4}$, $r^0_J(t)=\sum_{n=0}^{J-1} r^0_n\phi_n(t)$, we get that
\[
\| C^{1/2} r^0_J \|_{H^{(1-\alpha)/2}(0,t)} \leq \frac14|C^{1/2}\dot{e}_1|\sum_{n=0}^J b_n^n \|\phi_n\|_{H^{(1-\alpha)/2}(0,T)}\,,
\]
where 
\[
\|\phi_n\|_{H^{(1-\alpha)/2}(\mathbb{R})} = \Bigl(\int_\mathbb{R} (1+\omega^2)^{(1-\alpha)/2} |\mathcal{F}\phi_n(\omega)|^2\, d\omega\Bigr)^{1/2}
\]
with 
\[
\begin{aligned}
&|\mathcal{F}\phi_n(\omega)|^2 
= \left|\frac{1}{\sqrt{2\pi}\imath\omega}\Bigl(e^{-\imath\omega t_{n+}}-e^{-\imath\omega t_{n-}}\Bigr)\right|^2\\
&=\frac{1}{2\pi\omega^2}\Bigl((\cos(\omega t_{n+})-\cos(\omega t_{n-})^2+(\sin(\omega t_{n+})-\sin(\omega t_{n-})^2\Bigr)
\leq\frac{1}{\pi}\begin{cases}|\Delta t^2| \mbox{ for }\omega\in [-1,1]\\ 2\omega^{-2}\mbox{ for }|\omega|>1
\end{cases}\,,
\end{aligned}
\]
hence 
\[
\|\phi_n\|_{H^{(1-\alpha)/2}(\mathbb{R})}=
\Bigl(\int_{[-1,1]}+\int_{\mathbb{R}\setminus[-1,1]} (1+\omega^2)^{(1-\alpha)/2} |\mathcal{F}\phi_n(\omega)|^2\, d\omega\Bigr)^{1/2} 
\leq \frac{1}{\sqrt{\pi}}\bigl( 2\Delta t^2+ 2^{(5-\alpha)/2} \alpha^{-1}\bigr)^{1/2}\,.
\]
Moreover, 
\[
\sum_{n=0}^{J-1} b_n^n =
\frac{\Delta t^{1-\alpha}}{\Gamma(3-\alpha)} \Bigl( 1 + 2^{2-\alpha} - 2
 + \sum_{n=2}^{J-1} ((n+1)^{2-\alpha} - 2n^{2-\alpha} +(n-1)^{2-\alpha}) \Bigr)
\]
where by Taylor's theorem with $\xi_{n+}\in[n,n+1]$, $\xi_{n-}\in[n-1,n]$,
\[
\begin{aligned}
&\frac{\Delta t^{1-\alpha}}{\Gamma(3-\alpha)} \sum_{n=2}^{J-1} ((n+1)^{2-\alpha} - 2n^{2-\alpha} +(n-1)^{2-\alpha})= \frac{\Delta t^{1-\alpha}}{\Gamma(1-\alpha)} \sum_{n=2}^{J-1} ((n+\xi_{n+})^{-\alpha}+(n+\xi_{n-})^{-\alpha})
\\
&\leq \frac{2}{\Gamma(1-\alpha)} \, \Delta t \sum_{j=1}^{J} (\Delta t j)^{-\alpha}
\leq \frac{2}{\Gamma(1-\alpha)}\, \int_0^T x^{-\alpha}\, dx = \frac{2}{\Gamma(2-\alpha)} T^{1-\alpha}\,.
\end{aligned}
\]
Altogether we get for the initial data terms
\begin{equation}\label{estinitterms}
\begin{aligned}
&|M^{1/2} \hat{\dot{e}}_0 |^2  + | K^{1/2} e_0|^2 +  | K^{1/2} \Delta t \tilde{\dot{e}}_0|^2 
+\| C^{1/2} r^0_J \|_{H^{(1-\alpha)/2}(0,t)}^2\\
&=\Bigl( \|M^{1/2}(M^*)^{-1}C\dot{V}^\alpha\|_{L^\infty(0,T)}
+\|K^{1/2}(M^*)^{-1}C V^\alpha\|_{L^\infty(0,T)}
+\|M^{1/2}(M^*)^{-1}C\dot{V}^\alpha\|_{L^\infty(0,T)}\Bigr) O(\Delta t^2)\\
&=: c_0(V^\alpha) O(\Delta t^2)
\end{aligned}
\end{equation}

Thus the stabilty estimate Theorem \ref{th:stbdiscr} together with \eqref{estKterm}, \eqref{estCterm}, \eqref{estfterm}, \eqref{estinitterms}, yields
\[
\begin{aligned}
&\max_{n\in\{1\,\ldots,N\}} |M^{1/2} \hat{\dot{e}}_n |^2 + \| C^{1/2} \tilde{\dot{e}}_J \|_{H^{-(1-\alpha)/2}(0,t)}^2 
+ |K^{1/2} \Delta t \sum_{n=0}^N \tilde{\dot{e}}_n |^2\\
&\leq c(\alpha,T) \Bigl( c_0(V^\alpha) + \|M^{-1/2}K\ddot{u}\|_{L^\infty(0,T)} + \|M^{-1/2}C\ddot{V}^\alpha\|_{L^\infty(0,T)} + \|M^{-1/2}\ddot{f}\|_{L^\infty(0,T)}\Bigr) \Delta t^2 
\end{aligned}
\]
with $c_0(V^\alpha)$ as in \eqref{estinitterms}.

This is first of all an estimate for the velocity error only. In order to extract an estimate of the displacement error $e_n=u(t_n)-u_n$ as well, we make use of the identities
\[
\Delta t \sum_{n=0}^N \tilde{\dot{e}}_n = \Delta t \sum_{n=0}^N \tilde{v}_{cont,n} - \Delta t \sum_{n=0}^N \tilde{v}_n 
\]
where $\tilde{v}_{cont,n}=\frac{v(t_{n+1}) +2v(t_n) + v(t_{n-1})}{4}$, $n\in\{1,\ldots,N\}$, 
$\tilde{v}_{cont,0} = \frac{v(t_1)+v(t_0)}{2}$ and $v(t)=\dot{u}(t)$, so that
\[
\Delta t \sum_{n=0}^N \tilde{v}_{cont,n} = \frac{u(t_{N+1})+u(t_N)}{2}-u(t_0)
\]
and, due to \eqref{vnpdnp}
\[
\Delta t \sum_{n=0}^N \tilde{v}_n = \frac{d_{N+1}+d_N}{2}-d_0 = \frac{u_{N+1}+u_N}{2}-u_0\,,
\]
as well as $u_0=u(t_0)$.
 
\begin{theorem}\label{th:globerr}
For the time discretization scheme \eqref{eqdisc} with coefficients $b_l^n=\hat{b}_l^n$ according to \eqref{coeffGalerkin} with $f_{n+1}=f(t_{n+1})$, there exists a constant $C(T)$ depending only on $T$ such that the global error defined by \eqref{en}, \eqref{enhat}  satisfies
\[
\begin{aligned}
&\max_{n\in\{1\,\ldots,N\}} |M^{1/2} \hat{\dot{e}}_n |^2 + \| C^{1/2} \tilde{\dot{e}}_J \|_{H^{-(1-\alpha)/2}(0,t)}^2 
+ |K^{1/2} \frac{e_N+e_{N+1}}{2} |^2\\
&\leq c(\alpha,T) \Bigl( c_0(V^\alpha) + \|M^{-1/2}K\ddot{u}\|_{L^\infty(0,T)} + \|M^{-1/2}C\ddot{V}^\alpha\|_{L^\infty(0,T)} + \|M^{-1/2}\ddot{f}\|_{L^\infty(0,T)}\Bigr) \Delta t^2 
\end{aligned}
\]
for some constant $c(\alpha,T)>0$ and $c_0(V^\alpha)>0$ as in \eqref{estinitterms} with $V^\alpha=\partial_t^\alpha u$, provided $f$ and the solution $u$ of \eqref{eqn:semidiscr} exhibit the regularity required for finiteness of $c_0(V^\alpha)$, $\|M^{-1/2}K\ddot{u}\|_{L^\infty(0,T)}$, $\|M^{-1/2}C\ddot{V}^\alpha\|_{L^\infty(0,T)}$, $\|M^{-1/2}\ddot{f}\|_{L^\infty(0,T)}$. 
\end{theorem}

\begin{remark}

We have based our error analysis on the classical stability estimate approach of testing the wave type equation with $\dot{u}$ 
and therewith obtained an error estimate in the (discrete) energy norm for the wave equation, that is, $W^{1,\infty}(0,T;L^2(\Omega))\cap L^\infty(0,T;H^1(\Omega))$.
For doing so, it was crucial to have a coercivity estimate on the discretized fractional damping term, which we obtained by using the Galerkin discretization of the Abel integral operator from \cite{vogeli:disabel}. 
We are aware of the fact that the assumed smoothness of the exact solution in Theorem \ref{th:globerr} might be achievable only with rather smooth initial data.
In order to arrive at an $O(\Delta t^2)$ rate in the discrete $L^\infty(0,T;L^2(\Omega)$ norm for low regularity solutions, numerical schemes based on convolution quadrature and analysis approaches working in Laplace domain \cite{jin:diffusionwave} after transferring the Caputo derivative to a Riemann-Liouville one, might be fruitful.
\end{remark}

\section{Adjoint model}\label{sec:adjoint} 
\def\press{u}
To carry out a gradient descent method for minimizing the cost function in \eqref{MAP}, according to \eqref{gradJ} we have to apply the adjoint $\mathcal{G}^*$ of the forward operator $\mathcal{G}$ in appropriate function spaces (cf. \eqref{forwardop} below). The latter is defined via the solution $u$ to the fractional order initial boundary value problem
\begin{equation}\label{model} 
\left. \begin{array}{rcl}
\press_{tt} -b \Delta \partial_t^{\alpha} \press -c^2 \Delta \press &= 0  \quad &\mbox{ in } D \times [0, \infty) \\
\press &= 0 \quad & \mbox{ on } \partial D \times [0, \infty ) \\
\press(x,0) &= \bar{\press}_0(x) \quad &\mbox{ in } D \\
\press_t(x,0) &= 0 \quad &\mbox{ in } D
\end{array} \right\}
\end{equation}
after further mapping into the observations
\begin{equation*}
\press(x,t) = y(x,t) \quad \mbox{ on } \Sigma \times (0,T); \quad \Sigma = \partial \Omega.
\end{equation*}
Let $\mathcal{L}$ be the solution operator defined from $H_0^1(\Omega)$ to $W^{2,\infty} (0,\infty; (H^2 \cap H_0^1)(D)^*) \cap W^{1,\infty} (0,\infty; L^2(\Omega))\cap L^{\infty}(0,\infty; H_0^1(D))$ which maps the initial value $\press_0$ to the solution $\press$ of \eqref{model} with
\begin{equation*}
\bar{\press}_0 = \left\{ \begin{array}{cc} \press_0 & \mbox{ in } \Omega \\ 0 & \mbox{ in } D \setminus \Omega \end{array} \right.
\end{equation*}
From \cite[Proposition 3.1]{kaltenbacherrundell:fracpat} we know that the forward operator 
\begin{equation}\label{forwardop}
\mathcal{G}: H_0^1(\Omega) \to L^2(0,T; L^2(\Sigma))\,, \quad
\mathcal{G} = \operatorname{tr}_{\Sigma \times [0,T)} \circ \mathcal{L}.
\end{equation}
is well-defined and bounded.
To derive the adjoint $\mathcal{G}^*$ of $\mathcal{G}$, we first of all compute the Banach space adjoint of the time fractional derivative operator $\partial_t^{\alpha}$.
\[
\begin{aligned}
&\Gamma ( 1- \alpha ) \int_0^T (\partial_t^{\alpha} v)(t) u(t) \dd{t} = \int_0^T \int_0^t (t-s)^{-\alpha} v'(s) \dd{s} u(t) \dd{t} = \int_0^T v'(s) \int_s^T (t-s)^{-\alpha} u(t) \dd{t} \dd{s} = \\
&= \int_0^T v'(s) \left[ - \int_s^T \frac{(t-s)^{1-\alpha}}{1- \alpha} u'(t) \dd{t} + \frac{(T-s)^{1-\alpha}}{1- \alpha} u(T) \right] \dd{s} = \\
&= - \int_0^T v(s) \left[ \int_s^T (t-s)^{-\alpha} u'(t) \dd{t} - (T-s)^{-\alpha} u(T) \right]\dd{s} + v(T) \cdot 0 - v(0) \left[ - \int_0^T \frac{t^{1-\alpha}}{1-\alpha} u'(t) \dd{t} + \frac{T^{1-\alpha}}{1-\alpha} u(T) \right] = \\
&= - \int_0^T v(s) \left[ \int_s^T (t-s)^{-\alpha} u'(t) \dd{t} - (T-s)^{-\alpha} u(T) \right]\dd{s} -v(0) \int_0^T t^{-\alpha} u(t) \dd{t}
\end{aligned}
\]
So we have for $u,v \in H^1(0,T)$ that
\[
\begin{aligned}
&\int_0^T (\partial_t^{\alpha} v)(t) u(t) \dd{t} \\
&= -\frac{1}{\Gamma(1-\alpha)} \left\{ \int_0^T v(s) \left[ \int_s^T (t-s)^{-\alpha} u'(t) \dd{t} - (T-s)^{-\alpha} u(T) \right] \dd{s} + v(0) \left[ \int_0^T t^{-\alpha} u(t) \dd{t} \right] \right\} = \\
&= - \int_0^T v(s) ( \widetilde{\partial_t^{\alpha}} u)(s) \dd{s} - v(0) \tilde{I}^{\alpha} u
\end{aligned}
\]
with $\widetilde{\partial_t^{\alpha}}$ and $\tilde{I}^{\alpha}$ given by
\[
\begin{aligned}
(\widetilde{\partial_t^{\alpha}}u)(s) &= \frac{1}{\Gamma(1-\alpha)} \left[ \int_s^T (t-s)^{-\alpha} u'(t) \dd{t} - (T-s)^{-\alpha}u(T) \right] \\
\tilde{I}^{\alpha} u &= \frac{1}{\Gamma (1-\alpha)} \int_0^T t^{-\alpha} u(t) \dd{t}.
\end{aligned}
\]
Now we can also state the adjoint model of \eqref{model} and include our findings for the adjoint operator of the fractional derivative. To do so, we test the model \eqref{model} by some function $z$ and integrate by parts to obtain
\begin{equation}\label{ibp}
\begin{aligned}
&\int_0^T \int_{A} (\press_{tt} -b \Delta \partial_t^{\alpha} \press - c^2 \Delta \press)z\dd{x} \dd{t} \\
&= \int_0^T \int_A \press(s) \left\{ z_{tt} + b \Delta \widetilde{\partial_t^{\alpha}} z - c^2 \Delta z \right\} (s) \dd{x} \dd{s} + \left[ \int_A (\press_t z - \press z_t) \dd{x} \right]_0^T \\
&\quad+ c^2 \int_0^T \int_{\partial A} (- \partial_{\nu} \press z + \press \partial_{\nu} z) \dd{S(x)} \dd{t} \\
&\quad+ b \Bigl\{ \int_A \press(0)( \Delta \tilde{I}^{\alpha} z) \dd{x} + \int_0^T \int_{\partial A} \nu ( \widetilde{\partial_t^{\alpha}} z \nabla \press - \press \nabla \widetilde{\partial_t^{\alpha}}z) \dd{S(x)} \dd{t}\\
&\qquad\qquad+ \int_{\partial A} \nu ( \tilde{I}^{\alpha} z \nabla \press(0) - \press(0) \nabla \tilde{I}^{\alpha} z) \dd{S(x)} \Bigr\} 
\end{aligned}
\end{equation}

We subtract \eqref{ibp} with $A = D \setminus \Omega$ from \eqref{ibp} with $A=\Omega$ where we have the following situation: 
\begin{align*}
&A \in \{\Omega, D \setminus \Omega\}, \quad \nu_{\Omega} = - \nu_{D \setminus \Omega}, \quad \press = \mathcal{L} \press_0, \quad \mathcal{G} \press_0 =  \text{tr}_{I_{\Sigma \times (0,T) }}\press, \quad \bar{\press}_0|_{D \setminus \Omega} = 0, \quad \bar{\press}_0 |_{\Omega} = \press_0, \\
& \bar{\press}_0 |_{\partial \Omega} = 0 , \quad  \press|_{\partial \Omega} = 0 , \quad [\press] = 0, \quad [\partial_{\nu} \press] = 0 \text{ on } \Sigma \times (0,T), \quad \press_t(0)=0 \text{ on } D
\end{align*}

Given data $w \in L^2(0,T; \Sigma)$ this leads us to defining
\begin{equation}
\label{adjmodel1}
\left. \begin{array}{rcc}
z_{tt} + b \Delta \widetilde{\partial_t^{\alpha}} z - c^2 \Delta z &= 0 \quad &\mbox{ in } (\Omega \cup (D \setminus \Omega )) \times (0,T) \\
z&= 0 \quad &\mbox{ on } \partial D \times (0,T) \\
\left[ \partial_{\nu} (c^2  z - b \widetilde{\partial_t^{\alpha}} z) \right]  &= w \quad &\mbox{ on } \Sigma \times (0,T) \\
\left[z\right] &= 0 \quad &\mbox{ on } \Sigma \times (0,T) \\
z(T) &= 0 \quad &\mbox{ in } \Omega \cup D \setminus \Omega \\
z_t(T)  &= 0 \quad &\mbox{ in } \Omega \cup D \setminus \Omega
\end{array} \right\}
\end{equation}
such that \eqref{ibp} reduces to
\begin{equation}
0 = \int_{\Omega} \press_0 z_t(0) \dd{x} + \int_0^T \int_{\Sigma} \press w \dd{S(x)} \dd{t} + b \int_{\Omega} \press_0 ( \Delta\tilde{I}^{\alpha}z) \dd{x}
\end{equation}
which is equivalent to
\begin{equation}
\langle \mathcal{G} \press_0, w \rangle_{L^2(L^2)} = \int_{\Omega} \press_0 \left[ -z_t(0) - b( \Delta  \tilde{I}^{\alpha}z) \right] \dd{x} = \langle \press_0, \mathcal{G}^* w \rangle_{H_0^1(\Omega)}
\end{equation}
where $\tilde{z} := \mathcal{G}^* w$ solves
\begin{align*}
- \Delta \tilde{z} &= -z_t(0) -b \Delta \tilde{I}^{\alpha}z \quad &\mbox{ in } \Omega \\
\tilde{z} &= 0 \quad &\mbox{ on } \partial \Omega
\end{align*}

We then make a timeflip and consider $\bar{z}(t) := z(T-t)$, which due to the identity
\[
\begin{aligned}
& (\widetilde{\partial_t^{\alpha}}z)(T-t) = 
\frac{1}{\Gamma(1- \alpha)} \left[\int_{T-t}^{T} ( \tau -T+t)^{-\alpha} z_t(\tau) \dd{\tau} -t^{-\alpha} z(T)\right]\\
&=- \frac{1}{\Gamma(1-\alpha)}\left[\int_0^t (t-r)^{-\alpha} \bar{z}_t(r) \dd{r}-t^{-\alpha} z(T)\right]\\ 
&= - (\partial_t^{\alpha} \bar{z})(t) +\frac{1}{\Gamma(1-\alpha)}t^{-\alpha} \bar{z}(0)
\end{aligned}
\]
solves
\begin{align}
\label{adjmodel} 
\left. \begin{array}{rlc}
\bar{z}_{tt} - b \Delta \partial_t^{\alpha} \bar{z} - c^2 \Delta \bar{z} &= 0 \quad &\mbox{ in } (\Omega \cup (D \setminus \Omega )) \times (0,T) \\
\bar{z}&= 0 \quad &\mbox{ on } \partial D \times (0,T) \\
\left[\partial_{\nu} (c^2 \bar{z} + b \partial_t^{\alpha} \bar{z}) \right]  &= w(T-t) \quad &\mbox{ on } \Sigma \times (0,T) \\
\left[\bar{z}\right] &= 0 \quad &\mbox{ on } \Sigma \times (0,T) \\
\bar{z}(0) &= 0 \quad &\mbox{ in } \Omega \cup D \setminus \Omega \\
\bar{z}_t(0)  &= 0 \quad &\mbox{ in } \Omega \cup D \setminus \Omega
\end{array} \right\}
\end{align}
and $\tilde{z} = \mathcal{G}^* w$ solves 
\begin{equation}\label{ellbvp}
\begin{aligned}
- \Delta \tilde{z} &= \bar{z}_t(T) - b \Delta \hat{I}^{\alpha}\bar{z} \quad \mbox{ in } \Omega \\
z &= 0 \quad \mbox{ on } \partial \Omega
\end{aligned}
\end{equation}
with $\hat{I}^{\alpha} \bar{z} = \tilde{I}^{\alpha} z = \frac{1}{\Gamma(1- \alpha)}\int_0^T (T-s)^{-\alpha} \bar{z}(s) \dd{s}$. 

Therewith, the gradient $J'(u_0)$ cf. \eqref{gradJ} of $J$ can be computed by carrying out the following steps.
\begin{enumerate}
\item solve forward problem \eqref{model} on $D\times(0,T)$ and set $w=\Gamma_{noi}^{-1}(\operatorname{tr}_{\Sigma \times [0,T)} \press -y^\delta)$;
\item solve time reversed adjoint problem \eqref{adjmodel} on $D\times(0,T)$;
\item solve elliptic boundary value problem \eqref{ellbvp} on $\Omega$ and set
$J'(u_0)= \tilde{z}+\Gamma_{pr}^{-1}(u_0-u_0^*)$.
\end{enumerate} 
Both tasks 1. and 2. can be carried out by means of the method derived and analyzed in Section \ref{sec:timestepping}.
Similarly, also Hessian-vector products can be computed, cf., e.g., \cite{aaoBayes,flath:inversebayes}.

\section{Numerical reconstructions}\label{sec:reconstructions} 
In this section the results of the reconstruction by means of the devised numerical method are visualized. We perform the reconstructions in two space dimensions as relevant in medical imaging. To make use of Bayesian inference as described in section~\ref{sec:bayes} we employ the Python based software FEniCS \cite{AlnaesBlechta2015a} and hIPPYlib \cite{VillaPetraGhattas19}. We use the model \eqref{mod_init} with chosen parameters $c^2 = 1$ and $b=0.1$ and different values of $\alpha$, as well as observations on the boundary of a circle $\Sigma$, which is fully contained in the domain $D$. 
\begin{equation}
\begin{aligned}
u_{tt} - \Delta u - 0.1 \Delta \partial_t^{\alpha} u &= 0 \mbox{ in } D \times (0,T) \\
u(0) = u_0, \quad u_t(0) &= 0 \mbox{ in } D, \\
u &= y \quad \mbox{ on } \partial \Sigma \times (0,T)
\end{aligned}
\end{equation}
For our experiment we simulate observations on the boundary of the circle $\Sigma$. These observations with added white noise of level $\delta=0.01$ are used for the reconstructions. To avoid an inverse crime, we perform the reconstructions on a different discretization than the one we had used for the generation of the observations. We carry out tests for three different examples, where in each the searched for initial condition represents a constant inclusion in an otherwise homogeneous domain. The discretization is the same for all experiments and consists of linear Lagrangian finite elements on a mesh with $10,067$ points. We assume the domain $\Omega$ to be a square with side length $2$ centered around the origin. The observation boundary $\Sigma$ is a circle with radius $0.8$ centered at the origin. 
We consider the time interval $[0,1]$ which is discretized at points $t_i = i \cdot \Delta t$ with $\Delta t = 0.2$. The space for the reconstruction of the initial data is discretized with linear Lagrangian finite elements as well, but on a coarser mesh with $6,426$ points.
The setting remains the same for all examples, the results for the different initial conditions can be found in the following subsections.
In all examples we use a prior of the form 
\begin{equation}
\Gamma_{pr} = \left( \gamma I - \rho \Delta \right)^{-2},
\end{equation} 
where $I$ denotes the identity operator and $\Delta$ the Laplace operator on the given space. This is in fact the BiLaplacianPrior already implemented in hIPPYlib. 
In the tests below, the parameters $\gamma$ and $\rho$ are chosen ``by hand'' to yield good results; a more sophisticated choice can, e.g., be carried out by the discrepancy principle \cite{LuPereverzev2002}.


\subsection{Example 1}
\begin{figure}
\centering
\includegraphics[width=51mm]{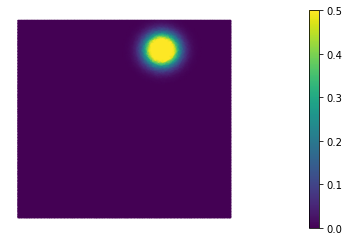}
\includegraphics[width=53mm]{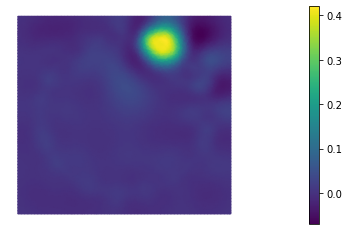}\\
\includegraphics[width=53mm]{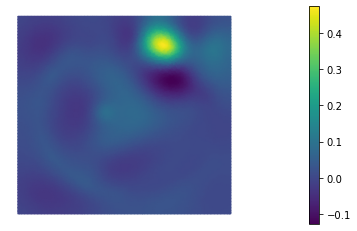}
\includegraphics[width=53mm]{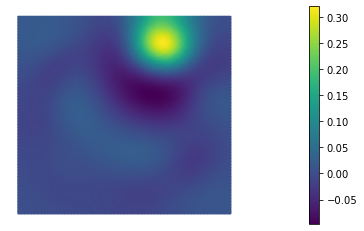}
\caption{Example 1: True initial condition (top left) and reconstructions with 
$\alpha=0.1$ (top right; $\gamma=10$, $\rho=0.03$), 
$\alpha=0.5$ (bottom left; $\gamma=10$, $\rho=0.03$)
$\alpha=0.9$ (bottom right; $\gamma=15$, $\rho=0.1$).
\label{fig_test1}}
\end{figure}

\begin{figure}
\centering
\includegraphics[width=0.7\textwidth]{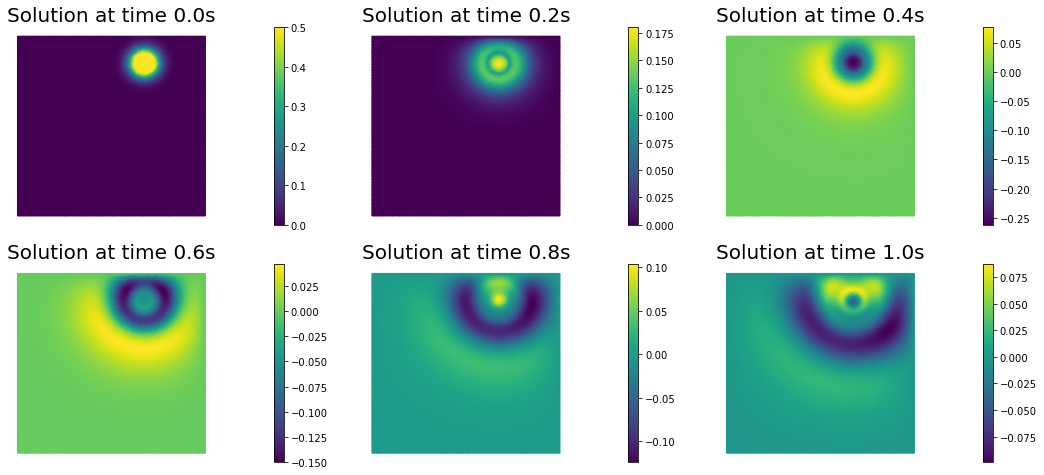}
\caption{Example 1: True state resulting from initial condition; $\alpha=0.5$.
\label{true_u1}}
\end{figure}

\begin{figure}
\centering
\includegraphics[width=0.7\textwidth]{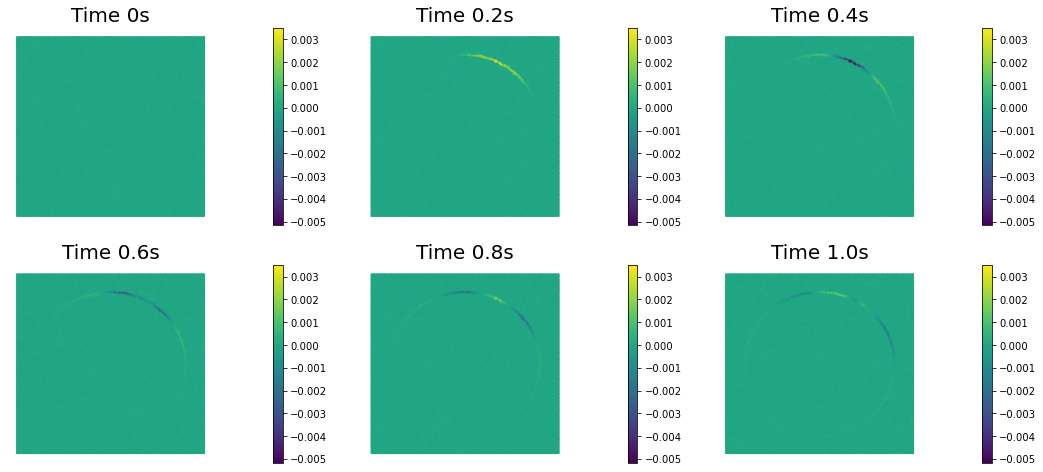}
\caption{Example 1: Observations on the boundary of $\Sigma$; $\alpha=0.5$.
\label{obs1}}
\end{figure}

For the first example we use an initial condition with an inclusion near the boundary of $\Sigma$, which, as one would expect, allows for the best reconstruction. 
The true initial condition and its reconstruction are given in Figure~\ref{fig_test1}. 
The fact that the problem is more ill-posed with stronger damping, that is, with larger $\alpha$, becomes evident from the reconstructions and from the fact that stronger regularization is needed.
Additionally the forward solution of the true initial condition and the resulting observations are given in Figures~\ref{true_u1} and~\ref{obs1}, respectively. Since these look quite similar for the other examples as well, we only include their visualization in this first example.

\subsection{Example 2}
\begin{figure}
\centering
\includegraphics[width=51mm]{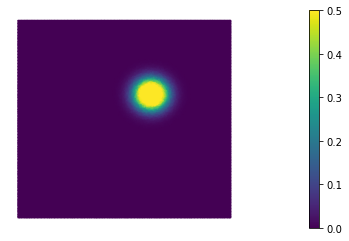}
\includegraphics[width=53mm]{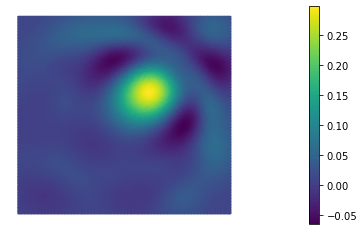}\\
\includegraphics[width=53mm]{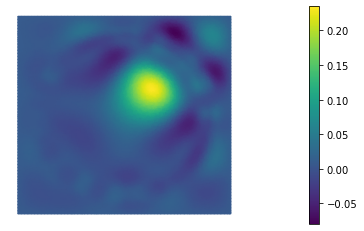}
\includegraphics[width=53mm]{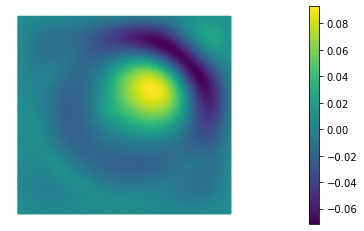}
\caption{Example 2: True initial condition (top left) and reconstructions with 
$\alpha=0.1$ (top right; $\gamma=10$, $\rho=0.01$), 
$\alpha=0.5$ (bottom left; $\gamma=15$, $\rho=0.01$)
$\alpha=0.9$ (bottom right; $\gamma=15$, $\rho=0.01$).
\label{fig:ex2}}
\end{figure}

For the second example we use an initial condition with an inclusion which is further away from the boundary $\Sigma$, which makes the reconstruction harder. 
The true initial condition and its reconstruction are given in Figure~\ref{fig:ex2}. It can be seen that the reconstruction gets more difficult the further away the inclusion is from the observation boundary. 

\subsection{Example 3}
\begin{figure}
\centering
\includegraphics[width=51mm]{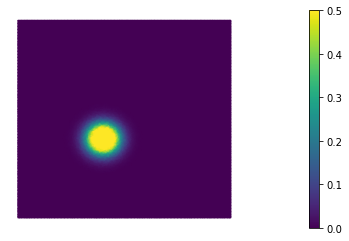}
\includegraphics[width=53mm]{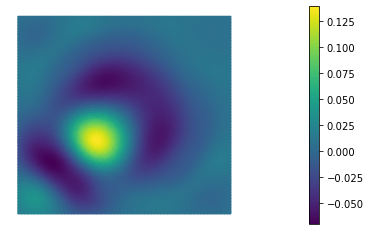}\\
\includegraphics[width=53mm]{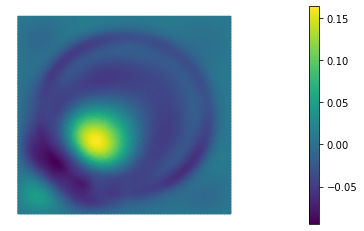}
\includegraphics[width=53mm]{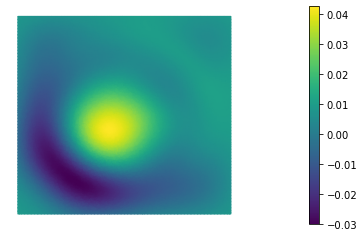}
\caption{Example 3: True initial condition (top left) and reconstructions with 
$\alpha=0.1$ (top right; $\gamma=15$, $\rho=0.1$), 
$\alpha=0.5$ (bottom left; $\gamma=15$, $\rho=0.1$) 
$\alpha=0.9$ (bottom right; $\gamma=15$, $\rho=0.1$).
\label{fig:ex3}}
\end{figure}

Finally, for the third and last example we use an initial condition with an inclusion which is almost in the center of the circle $\Sigma$. 
The true initial condition and its reconstruction are given in Figure~\ref{fig:ex3}. Here, the reconstruction is worse than  in the examples above (also with respect to the actual value of the inclusion); additionally, we can clearly see the observation circle as an image artefact.

\section{Conclusions and Outlook}
In this paper we have discussed two tasks that are crucial for the numerical solution of the inverse problem of photoacoustic or thermoacoustic tomography in the presence of fractional derivative attenuation. For the forward problem of solving a second order wave equation with fractional damping, we have derived a time stepping scheme based on the Newmark scheme that is modified in order to take the time fractional derivative term into account.
Using the Galerkin scheme from \cite{vogeli:disabel} for the Abel integral operator contained in the fractional derivative, we obtain a stability result for the time discretization, as well as a convergence rate with respect to the discrete energy norm. 
Our second key contribution for the efficient solution of the inverse problem consists of the derivation of an adjoint scheme for gradient computation, which is employed within a minimization based regularization method. Also here, the fractional derivative term requires extra treatment.

Future work will be concerned with numerical solution methods for this inverse problem in time domain in the context of  alternative damping models, cf., e.g., \cite{kaltenbacherrundell:fracpat}, see also \cite{BakerBanjai20}. For some frequency domain based approaches we refer to, e.g., \cite{ElbauScherzerShi:2017,KowarScherzer:2011}.

\section*{Appendix}
\subsection*{Proof of Lemma~\ref{lemma1}}
We have 
\begin{equation}
\label{eqpr1}
\left(w_{n+1} - w_{n-1}, \frac{w_{n+1} + 2 w_n + w_{n-1}}{4}\right) = \frac{1}{4} \left[ |w_{n+1}|^2 - |w_{n-1}|^2 + 2(w_{n+1}, w_n) - 2(w_{n-1}, w_n) \right]
\end{equation}
where it holds $2(w_{n+1}, w_n) = -|w_{n+1} - w_n|^2 + |w_{n+1}|^2 + |w_n|^2$ and $-2(w_{n-1}, w_n) = |w_{n-1} - w_n|^2 - |w_{n-1}|^2 - |w_n|^2$, so that \eqref{eqpr1} is equal to
\[ \begin{aligned}
&\frac{1}{2} \left( |w_{n+1}|^2 - |w_{n-1}|^2 \right) - \frac{1}{4} \left( |w_{n+1} - w_n |^2 - |w_{n-1} -w_n |^2 \right) =\\
& \frac{1}{2} \left( |w_{n+1}|^2 + |w_n|^2 -(|w_n|^2+ |w_{n-1}|^2) \right) - \frac{1}{4} \left( |w_{n+1} - w_n |^2 - |w_{n-1} -w_n |^2 \right) 
\end{aligned} \]
Now we can compute the sum to obtain
\[ \begin{aligned}
&\sum_{n=1}^N \left(w_{n+1} - w_{n-1}, \frac{w_{n+1} + 2w_n + w_{n-1}}{4}\right) =  \\
&\frac{1}{2}( |w_{N+1}|^2 + |w_N|^2 -|w_1|^2 -  |w_0|^2) + \frac{1}{4} ( |w_1 - w_0|^2 - |w_{N+1} - w_N|^2) = \\
&= \frac{1}{4} |w_{N+1}|^2 + \frac{1}{4}|w_N|^2 + \frac{1}{2}(w_{N+1}, w_N) - \left[ \frac{1}{4}|w_1|^2 + \frac{1}{4}|w_0|^2 + \frac{1}{2} (w_1, w_0) \right]= \\
& = \left|\frac{w_{N+1} + w_N}{2} \right|^2 - \left|\frac{w_1 + w_0}{2}\right|^2.
\end{aligned} \]

\subsection*{Proof of Lemma~\ref{lemma2}}
We have
\begin{equation}
\begin{aligned}
\label{eqpr2}
\sum_{n=0}^N \left(\sum_{j=1}^n \frac{\tilde{w}_j + \tilde{w}_{j-1}}{2}, \tilde{w}_n \right) =\frac{1}{2} \left[ \sum_{n=0}^N \sum_{j=1}^n (\tilde{w}_j, \tilde{w}_n ) + \sum_{n=0}^N \sum_{i=0}^{n-1} ( \tilde{w}_i, \tilde{w}_n) \right]\,,
\end{aligned}
\end{equation}
where the first sum can be written as
\begin{equation*}
\sum_{n=0}^N \sum_{j=1}^n (\tilde{w}_j, \tilde{w}_n ) = \sum_{n=0}^N | \tilde{w}_n |^2 + \sum_{n=0}^N \sum_{j=1}^{n-1} (\tilde{w}_j, \tilde{w}_n) = \sum_{n=0}^N |\tilde{w}_0|^2 + \sum_{n=0}^N \sum_{j=0}^{n-1} ( \tilde{w}_j, \tilde{w}_n) - (\tilde{w}_0, \sum_{n=0}^N \tilde{w}_n)
\end{equation*}
and the second sum is reformulated as
\begin{equation*}
 \sum_{n=0}^N \sum_{i=0}^{n-1} ( \tilde{w}_i, \tilde{w}_n) = \frac{1}{2} \sum_{n=0}^N \sum_{\substack{i=0 \\ i \neq n} }^N (\tilde{w}_i, \tilde{w}_n )
\end{equation*}
So in total we have that (\ref{eqpr2}) is equal to
\[
\begin{aligned}
&\frac{1}{2} \left\{\sum_{n=0}^N | \tilde{w}_n |^2 + \sum_{n=0}^N \sum_{\substack{i=0 \\ i \neq n}}^N (\tilde{w}_i, \tilde{w}_n ) - (\tilde{w}_0, \sum_{n=0}^N \tilde{w}_n) \right\} = \frac{1}{2} \left\{ \sum_{n=0}^N \sum_{i=0}^N (\tilde{w}_i, \tilde{w}_n) - (\tilde{w}_0, \sum_{n=0}^N \tilde{w}_n)\right\} = \\
& = \frac{1}{2} \left\{ \left| \sum_{n=0}^N \tilde{w}_n \right|^2 - (\tilde{w}_0, \sum_{n=0}^N\tilde{w}_n) \right\} \geq \frac{1}{4} \left\{ \left| \sum_{n=0}^N \tilde{w}_n \right|^2 - | \tilde{w}_0|^2 \right\}
\end{aligned}
\]
where we used Young's inequality as last step.


\end{document}